\documentclass[a4paper]{article}
\usepackage{amssymb,amsbsy,amsmath,amsfonts,amssymb,amscd}
\usepackage{latexsym}

\usepackage{latexsym}

\newcommand{\lie}[1] {\mathfrak{#1}}  
\newcommand{\bb}[1]{{\mathbb #1}}    

\newcommand{\bbR}{{\bb R}}

\newcommand{\bbZ}{{\bb Z}} 
\newcommand{\bbQ}{{\bb Q}}

\newcommand{\bbC}{{\bb C}}

\newcommand{\rmod}{ / \,}

\newcommand{\GL}{{\rm GL}}

\newcommand{\Aff}{{\rm Aff}}
\newcommand{\Aut}{{\rm Aut}}

\newcommand{\Z}{{\rm Z}} 
\newcommand{\Fitt}{{\rm Fitt}} 
\newcommand{\ur}{{\rm u}} 
\newcommand{\C}{{\rm Z}} 

\newcommand{\bF}{{\mathbf{F}}}

\newcommand{\bU}{{\mathbf{U}}}
\newcommand{\bG}{{\mathbf{G}}}
\newcommand{\bH}{{\mathbf{H}}}

\newcommand{\bN}{{\mathbf{N}}}

\newcommand{\bD}{{\mathbf{D}}}
\newcommand{\bT}{{\mathbf{T}}}
\newcommand{\bI}{{\mathbf{I}}}

\newcommand{\bX}{{\mathbf{X}}}

\newcommand{\lu}{{\lie{u}}}
\renewcommand{\lg}{{\lie{g}}}

\newcommand{\rank}{{\rm rank }\, }

{ \par  \medskip  \noindent  {\bf Definition} \hspace{0.5em} }%
{\par \medskip } 

\newtheorem{example1}{Example}[section]

{ \begin{example1} \rm }%
{ \end{example1}} 
\newenvironment {proof}%
{ \noindent {\em Proof. }}%
{\hspace*{\fill}$\Box$\par \medskip } 
\newenvironment{prf}[1]%
{ \noindent {\em #1 \hspace{0.5em}} }%
{\hspace*{\fill}$\Box$\par \medskip } 
\newenvironment {remark}%
{{\em Remark \hspace{0.5em}} }%
{\par \medskip }

\newtheorem{proposition}{Proposition}[section]
\newtheorem{definition1}[proposition]{Definition}
\newtheorem{theorem}[proposition]{Theorem}
\newtheorem{lemma}[proposition]{Lemma}
\newtheorem{corollary}[proposition]{Corollary}
\newenvironment{definition}%
{ \begin{definition1} \rm }%
{ \end{definition1}}

\newcommand{\ac}[1]{\overline{#1}} 
\renewcommand{\rank}{{\rm rank}\,}
\newcommand{\hol}{{hol}}

\date{November 12, 2003}

\title{Infra-Solvmanifolds and Rigidity of Subgroups
in Solvable Linear Algebraic Groups}

\author{Oliver Baues \thanks{e-mail: oliver@math.ethz.ch} \\
Departement Mathematik\\ 
ETH-Zentrum \\ 
R\"amistrasse 101 \\
CH-8092 Z\" urich}

\begin{document}
\maketitle 
\begin{abstract} 
We give a new proof that compact infra-solvmanifolds 
with isomorphic fundamental groups are
smoothly diffeomorphic. More generally, we
prove rigidity results for manifolds which are
constructed using affine actions of virtually polycyclic
groups on solvable Lie groups. Our results are derived from
rigidity properties of subgroups in solvable
linear algebraic groups. 
\end{abstract}

\section{Introduction}
A closed manifold $M$ is called 
topologically rigid if every homotopy equivalence 
$h\colon N \to M$ 
from another manifold $N$ is homotopic to a homeomorphism. 
The Borel conjecture 
expects every closed aspherical 
manifold to be topologically rigid. 
The manifold $M$ is called smoothly  
rigid if every homotopy equivalence is
homotopic to a diffeomorphism. 
Geometric methods are useful to prove smooth 
rigidity inside some classes of closed 
aspherical manifolds. Well known cases are, for example,  
locally symmetric spaces of non-compact type \cite{Mostow4}, 
or flat Riemannian manifolds \cite{Bieberbach}. 
In this paper, we study the smooth rigidity 
problem for infra-solvmanifolds. These manifolds 
are constructed by considering isometric affine actions
on solvable Lie groups.
  
The fundamental group of an infra-solvmanifold
is a virtually polycyclic group.
A result of Farrell and Jones \cite{FJ1}
on aspherical manifolds with virtually 
polycyclic fundamental group shows 
that infra-solvmanifolds are topologically rigid. 
Yet, an argument due to Browder \cite{Browder} implies
that there exist smooth manifolds which are homeomorphic 
but not diffeomorphic to the $n$-torus, for
$n \geq 5$. \mbox{Farrell} and Jones \cite{FJ2} proved that  
any two compact infrasolvmanifolds  
of dimension not equal to four, 
whose fundamental groups are isomorphic, are diffeomorphic.
This generalizes previous results of 
Bieberbach \cite{Bieberbach} on compact flat Riemannian manifolds, 
Mostow \cite{Mostow1} on compact solv-manifolds 
and of Lee and Raymond \cite{Lee-Raymond} on infra-nilmanifolds.  
The proof of Farrell and Jones requires smoothing theory and
the topological rigidity result. Recent results of 
Wilking on rigidity properties of isometric actions
on solvable Lie-groups \cite{Wilking} imply the
smooth rigidity of infra-solvmanifolds in all
dimensions, giving an essentially geometric
proof. 

In well known cases, smooth rigidity properties of 
geometric manifolds are closely 
connected to rigidity properties
of lattices in Lie groups.
The aim of the present paper is to establish the 
smooth rigidity of infra-solvmanifolds 
from natural rigidity properties of 
virtually polycyclic groups in linear 
algebraic groups. More generally, we prove rigidity 
results for manifolds 
which are constructed using affine, not necessarily isometric, 
actions of virtually polycyclic groups on solvable Lie groups. 
This approach leads us
to a new proof of the rigidity of infra-solvmanifolds, 
and also to a geometric characterization of infra-solvmanifolds
in terms of polynomial actions on affine space
$\bbR^n$. As an application of the latter 
point of view we compute the cohomology
of an infra-solvmanifold using the 
finite-dimensional complex of polynomial differential forms. 
This generalizes a result of Goldman \cite{Goldman} on 
compact complete affine manifolds. As another application, 
we show that every infra-solvmanifold has maximal torus-rank. 
Our approach towards rigidity of infra-solvmanifolds
also suggests to study the rigidity-problem for the potentially 
bigger class of manifolds which are constructed using affine  
actions of virtually polycyclic groups on solvable 
Lie groups. Our main result establishes smooth rigidity for 
virtually polycyclic affine actions if the
holonomy of the action is contained in a reductive group,
generalizing the particular case of isometric actions. 

\paragraph{Infra-solvmanifolds}
We come now to the definition of infra-solvmanifolds.
Let $G$ be a Lie-group and let $\Aff(G)$ denote the 
semi-direct product $G \rtimes \Aut(G)$, where $\Aut(G)$
is the group of automorphisms of $G$. We view $\Aff(G)$
as a group of transformations acting on  $G$. If $\Delta$ is a
subgroup of $\Aff(G)$ then let $\Delta_0$ denote its 
connected component of identity, and 
$\hol(\Delta) \leq \Aut(G)$ its image under
the natural homomorphism $\Aff(G) \rightarrow \Aut(G)$.  

\begin{definition} \label{infrasolv}
An {\em infra-solvmanifold\/} is a manifold of the form 
$\Delta \backslash G$, where $G$ is a connected, 
simply connected solvable Lie group, 
and $\Delta$ is a torsion-free subgroup of ${\rm Aff}(G)$ which satisfies
(1) the closure of $\hol(\Delta)$ in ${\rm Aut}(G)$ is compact.
\end{definition} 


The manifold $\Delta\backslash G$ is a smooth manifold 
with universal cover diffeomorphic to $\bbR^m$, 
$m=\dim G - \dim \Delta_0$, where $\Delta_0$ is
the connected component of identity in $\Delta$. 
The fundamental group 
of $\Delta \backslash G$ is isomorphic to 
$\Gamma= \Delta /\Delta_0$. 
It is known that $\Delta\backslash G$
is finitely covered by a solv-manifold, i.e., a 
homogeneous space of a solvable Lie group.
By a result of Mostow \cite{Mostow3} 
the torsion-free group $\Gamma$ is 
then a virtually polycyclic group. 
(Recall that a group $\Gamma$ is
called virtually polycyclic (or polycyclic by finite)
if it contains a subgroup 
$\Gamma_0$ of finite index which is polycylic, i.e., 
$\Gamma_0$ admits a finite normal series with
cyclic quotients. The number of infinite cyclic factors
in the series is an invariant of $\Gamma$ called the
rank of $\Gamma$.) If $\Delta\backslash G$
is compact then $\dim \Delta\backslash G$ 
equals the rank of $\Gamma$. 
Not every smooth manifold which is finitely covered by
a compact solvmanifold is diffeomorphic to an 
infra-solvmanifold. By the work of Wall 
et al.\ (see \cite{KS}), there exist fake tori which are 
finitely covered by standard tori. 
Consequently, these smooth manifolds do not 
carry any infra-solv structure. 
 
\paragraph{Main results}
Let $\Gamma$ be a torsion-free virtually 
polycyclic group. To $\Gamma$ we associate in a functorial way a
solvable by finite real linear algebraic group $H_\Gamma$ 
which contains $\Gamma$ as a discrete and Zariski-dense subgroup. 
The group $H_\Gamma$ is called the real algebraic hull 
for $\Gamma$. The construction of the algebraic hull 
for $\Gamma$ extends results of Malcev \cite{Malcev} 
on torsion-free nilpotent groups, 
and results of Mostow \cite{Mostow2} on  
torsion-free polycyclic groups. 
The extended construction was 
first announced in \cite{Baues1}. 
The details are provided in Appendix A 
of this paper.  

We explain now the role the real algebraic hull 
plays in the construction of infra-solvmanifolds.
Let $T \leq H_\Gamma$ be a maximal reductive subgroup,
and let $U$ denote the unipotent radical of $H_\Gamma$. 
Then $H_\Gamma$ decomposes as a semi-direct product 
$H_\Gamma = U \cdot T$. 
The splitting induces an injective  homomorphism 
$\alpha_T: H_\Gamma \rightarrow \Aff(U)$
and a corresponding affine action of 
$\Gamma \leq H_\Gamma$ on $U$.  
The quotient space 
$$ M_\Gamma  \; = \; {_{\displaystyle \alpha_T(\Gamma)}}  \backslash  U$$
is a compact aspherical manifold of dimension $n= \rank \Gamma$,
and has universal cover $U = \bbR^n$. In fact, 
we show that $M_\Gamma$ is an infra-solvmanifold. 
We call every manifold $M_\Gamma$ which arises by 
this construction a standard $\Gamma$-manifold.

We prove: 

\begin{theorem} \label{main1} 
Let\/ $\Gamma$ be a torsion-free virtually 
polycyclic group. Then $M_\Gamma$ is a
compact infra-solvmanifold
and the fundamental group $\pi_1(M_\Gamma)$ is 
isomorphic to $\Gamma$.
Every two standard $\Gamma$-manifolds 
are diffeomorphic and 
every given isomorphism of fundamental 
groups of standard $\Gamma$-manifolds is 
induced by a smooth diffeomorphism.
\end{theorem}

Let $G$ be a connected, simply connected 
Lie group, and let $\lie{g}$ denote its
Lie algebra.  The group $\Aut(G)$ attains the structure
of a real linear algebraic group since it has
a natural identification with the group $\Aut(\lie{g})$
of Lie algebra automorphisms of $\lie{g}$. Our main
result is: 

\begin{theorem} \label{main} 
Let $G$ be a connected, simply connected solvable 
Lie group. Let $\Delta \leq \Aff(G)$ be a solvable by finite 
subgroup which acts freely and properly on $G$ with compact
quotient manifold $M = \Delta \backslash G$. 
Assume that one of the following two conditions is 
satisfied:
\begin{itemize}
\item[i)] $G$ is nilpotent, or 
\item[ii)] $hol(\Delta) \leq \Aut(G)$ 
is contained in a reductive subgroup of $\Aut(G)$. 
\end{itemize}
Then\/ the group $\Gamma = \Delta/\Delta_0$
is virtually polycyclic, and $M$ is diffeomorphic 
to a standard $\Gamma$-manifold. 
\end{theorem}  

We deduce: 

\begin{theorem} \label{main2} 
Every compact infra-solvmanifold is 
smoothly diffeomorphic to a standard $\Gamma$-manifold.  
\end{theorem} 

\begin{corollary} \label{srigid}
Compact infra-solvmanifolds are 
smoothly rigid. In particular, 
every two compact infra-solvmanifolds with
isomorphic fundamental groups are smoothly diffeomorphic.  
\end{corollary}

Theorem \ref{main1} also implies
the following result which was first proved by Auslander 
and Johnson \cite{AJ}. Their construction is
different from ours. 
  
\begin{corollary} Every torsion-free virtually polycyclic group
is the fundamental group of a compact infra-solvmanifold.
\end{corollary}


The torus rank $r$ of a manifold $M$ is
the maximum dimension of a torus 
which acts almost freely and 
smoothly on $M$. 
For a closed aspherical manifold $M$,
$r$ is bounded by the rank
of the center of the fundamental group.
If $r$ equals  the rank of the center then the 
torus rank of $M$ is said to be maximal.  
It is known (see \cite{Lee-Raymond2}) that 
that the torus rank of a solvmanifold is maximal
and the result is expected to hold for 
infra-solvmanifolds as well. It is 
straightforward to see that standard 
$\Gamma$-manifolds admit maximal torus actions.
Therefore, we also have: 

\begin{corollary} \label{torusrank}
Every infra-solvmanifold
has maximal torus rank.
\end{corollary}

Let $U$ be a connected, simply connected, nilpotent
Lie group, and let $\lu$ denote its Lie-algebra. 
Nomizu \cite{Nomizu} proved that the cohomology of a
compact nilmanifold $M= U/\Gamma$, where 
$\Gamma \leq U$ is a lattice, is isomorphic 
to the cohomology of the complex of 
left invariant differential forms on $U$.  
This means that the cohomology of 
the nilmanifold $M$ is computed by the Lie algebra 
cohomology $H^*(\lu)$.
Now let $\Gamma$ be a torsion-free virtually 
polycyclic group and $M_\Gamma$ a standard $\Gamma$-manifold.
Let $H_\Gamma=U \cdot T$ be the real algebraic hull
for $\Gamma$, where $U$ is the unipotent radical and
$T$ is maximal reductive. Then $T$ acts by automorphisms on $U$ 
and on the cohomology  ring $H^*(\lu)$.
Let $H^*(\lu)^T$ denote the $T$-invariants in $H^*(\lu)$.  
Let $M$ be an infra-solvmanifold 
with fundamental group $\Gamma$.
By Theorem \ref{main2}, $M$ is 
diffeomorphic to the standard $\Gamma$-manifold $M_\Gamma$.
Hence, the following 
result computes the cohomology of $M$: 

\begin{theorem} \label{cohomology1}
Let $M_\Gamma$ be a standard\/ $\Gamma$-manifold.
Then the de Rham-co\-ho\-mo\-lo\-gy 
ring $H^*(M_\Gamma)$ is isomorphic 
to $H^*(\lu)^T$. 
\end{theorem} 

We remark that the theorem implies
that the discrete group cohomology 
of $\Gamma$, 
$H^*(\Gamma, \bbR) = H^*(M_\Gamma)$,  is 
isomorphic to the rational cohomology 
(see \cite{Hochschild}[Theorem 5.2])
of the real linear algebraic 
group $H_\Gamma$.

\paragraph{Some historical remarks}
We want to give a few more 
historical remarks about the context 
of our paper, and the techniques we use. As our main
tool we employ the algebraic hull functor 
which naturally associates a linear algebraic
group to a (torsion-free) virtually polycyclic group or 
to a solvable Lie group. This functor was 
considered by Mostow in his paper \cite{Mostow2}. 
Auslander and Tolimieri solved the main open problems
on solv-manifolds at their time 
using the technique of the nilpotent shadow 
and semi-simple splitting for solvable Lie groups
(see \cite{AusTol,Auslander}). 
Mostow remarked then in \cite{Mostow5}
that the nilpotent shadow and 
splitting construction may be derived 
naturally from the algebraic hull, and 
reproved the Auslander-Tolimieri results,
as well as his older result on the 
rigidity of compact solv-manifolds. 
In our paper, we establish and use the properties of 
the algebraic hull functor for the class of 
virtually polycyclic groups not containing finite 
normal subgroups. We provide the necessary 
results and proofs about the hull functor 
in an appendix. Immediate applications are 
then our rigidity results and cohomology 
computations for infra-solvmanifolds. 
In \cite{Baues1} we give another
application of the hull functor in the context of 
affine crystallographic groups
and their deformation spaces.   

\paragraph{Arrangement of the paper}
We start in \S 2 with some preliminaries 
on real algebraic and syndetic hulls for virtually 
polycylic groups, affine actions and splittings
of real algebraic groups. The necessary results
about the construction of algebraic hulls are provided in
Appendix A. In \S 3 we prove Theorem \ref{main1} and
Theorem \ref{main}. In \S 4 we provide some applications
on the geometry of infrasolvmanifolds. 
In particular, 
we show that infra-solvmanifolds are distinguished 
in the class of aspherical compact differentiable 
manifolds with a virtually polycyclic 
fundamental group by the existence of
a certain atlas whose coordinate changes 
are polynomial maps. As an application, 
we compute the cohomology of infra-solvmanifolds
in terms of polynomial differential forms,  
and derive Theorem \ref{cohomology1}. 

\paragraph{Acknowledgement}
I thank Fritz Grunewald, Burkhard Wilking and 
Wilhem Singhof for helpful comments on an
earlier draft of this article.

\section{Hulls and splittings}
We need some terminology concerning real algebraic 
groups. For terminology on algebraic groups see
also Appendix \ref{appendixA}. 
Let $\bG$ be a $\bbR$-defined linear 
algebraic group. The 
group of real points $G= \bG_\bbR \leq \GL_n(\bbR)$ will 
be called a real algebraic group.
The group $G$ has the natural Euclidean topology which 
turns it into a real Lie-group but it carries 
also the Zariski-topology induced from $\bG$.
Let $H= \bH_\bbR$ be another real algebraic group. 
A group homomorphism 
$\phi: G \rightarrow H$ 
is called an algebraic homomorphism
if it is the restriction of a $\bbR$-defined morphism 
$\bG \rightarrow \bH$ of linear algebraic groups. 
If $\phi$ is an isomorphism of groups which is algebraic 
with algebraic inverse, then $\phi$ is called an 
algebraic isomorphism.  
We let $G^0$ denote the Zariski-irreducible
component of identity in $G$, and $G_0$ the connected component 
in the Euclidean topology. 
In particular, $G_0 \leq G^0$ is
a subgroup of finite index in $G$.  
If $g$ is an element of $G$ then 
$g = g_u g_s$ denotes the Jordan-decomposition
of $g$. Here $g_u \in G$ is unipotent, $g_s \in G$ is
semisimple, and $g_u ,g_s$ commute. 
Let $M \subset G$ be a subset. 
Then $\ac{M}$ denotes the Zariski-closure
of $M$ in $G$.
We put $M_u = \{ g_u \mid g \in M \}$,
$M_s=  \{ g_s \mid g \in M \}$. 
We let $\ur(G)$ denote the unipotent 
radical of $G$, i.e., the maximal normal subgroup of $G$ which
consists of unipotent elements. 
 
\subsection{Solvable by finite real algebraic groups}
A linear algebraic group $\bH$ is called   
solvable by finite if $\bH^0$ is solvable.
Assume that $\bH$ is 
solvable by finite. Then $\bH_u = \ur(\bH)$. 
In particular, for any subgroup $G$ of $\bH$, 
$\ur(G) = G \cap G_u$. If $G$ is a nilpotent
subgroup then (compare, \cite[\S 10]{Borel}) 
$G_u$ and $G_s$ are subgroups of $\bH$, and 
$G \leq G_u \times G_s$. A Zariski-closed 
subgroup $\bT \leq \bH $ which consists only of 
semi-simple elements is called a $d$-subgroup
of $\bH$. The group $H= \bH_\bbR$  is called a solvable 
by finite real algebraic group. Every Zariski-closed 
subgroup $T \leq H$ consisting of semi-simple 
elements is called a $d$-subgroup of $H$. Any $d$-subgroup 
of $H$ is an abelian by finite group, and
its identity component $T^0$ is a
real algebraic torus. 
   
\begin{proposition} \label{Usplitting} 
Let $H$ be a solvable 
by finite real linear algebraic group. 
Let $T$ be a maximal $d$-subgroup of $H$,
and $U= \ur(H)$ the unipotent 
radical of $H$. Then 
\begin{equation*}  
H = U \cdot T \; \text{ (semi-direct product)} \; . 
\end{equation*} 
Moreover, any two maximal $d$-subgroups $T$ and
$T'$ of $H$ are conjugate by an element of $U$.   
\end{proposition} 
\begin{proof} Let us assume that $H \leq \bH$
is a Zariski-dense subgroup.  
Let $\bT$ be the Zariski-closure
of $T$ in $\bH$. Then $\bT$ is a $\bbR$-defined 
subgroup of $\bH$, and a $d$-subgroup. 
Also $T \leq \bT_\bbR$ and, by maximality
of $T$, $T = \bT_\bbR$. Moreover, $\bT$ is a 
maximal reductive $\bbR$-defined subgroup of $\bH$. 
Therefore, by a
well known result (see \cite[Proposition 5.1]{BS}) 
$\bH = \bU \cdot \bT$, where $\bU = \ur(\bH)$, 
and every two $\bbR$-defined
maximal reductive subgroups $\bT$ and $\bT'$ are
conjugate by an element of $\bU_\bbR=U$.
Then the decomposition of $H$ 
follows. Since
$T$ and $T'$ are the group of real-points in 
maximal $d$-subgroups  $\bT$ and $\bT'$
they are conjugate by an element of $U$.  
\end{proof}

\subsection{Algebraic hulls} \label{sah}
Let $\Gamma$ be a torsion-free virtually polycyclic group.
We introduce the concept of an algebraic hull for $\Gamma$.
For more details and proofs see Appendix \ref{appendixA}.
Let $\bG$ be a linear algebraic group, and let 
$\bU$ denote the unipotent radical 
of $\bG$. We say that $\bG$ has a {\em strong unipotent 
radical\/} if the centralizer $\Z_{\bG}(\bU)$ is contained 
in $\bU$.  

\begin{theorem} \label{alghull0} 
There exists a $\bbQ$-defined linear algebraic group $\bH$ and
an injective homomorphism $\psi: \Gamma \longrightarrow \bH_{\bbQ}$ so that,
\begin{itemize}
\item[i)] $\psi(\Gamma)$ is Zariski-dense in $\bH$, 
\item[ii)] $\bH$ has a strong unipotent radical $\bU$,
\item[iii)] $\dim \bU = \rank \Gamma$. 
\end{itemize} 
\end{theorem}

We call the $\bbQ$-defined linear algebraic group $\bH$ 
the {\em algebraic hull for $\Gamma$}. The homomorphism  
$\psi$ may be chosen so that $\psi(\Gamma) \cap  \bH_{\bbZ}$ 
has finite index in  $\psi(\Gamma)$.  
Let $k \leq \bbC$
be a subfield. The hull $\bH$ together with a 
Zariski-dense embedding 
$\psi: \Gamma \longrightarrow \bH_k$ of $\Gamma$ into 
the group of $k$-points of $\bH$
satisfies  the following rigidity property:
\begin{itemize}
\item[(*)] 
Let $\bH'$ be another linear algebraic group and
$\psi': \Gamma \longrightarrow \bH'_k$ an injective
homomorphism so that i) to iii) above are satisfied with
respect to $\bH'$. Then there  
exists a $k$-defined isomorphism 
$\Phi: \bH \rightarrow \bH'$ so that 
$\psi' = \Phi \circ \psi$. 
\end{itemize}
In particular, 
the group $\bH$ is determined by the conditions
i)-iii) up to $\bbQ$-defined 
isomorphism of linear algebraic groups.   

\paragraph{The real algebraic hull for $\Gamma$}
Let $\bH$ be an algebraic hull for $\Gamma$,
$H = \bH_\bbR$ the group of real points. 
Put $U= \ur(H)$. 
Then there exists an injective homomorphism 
$\psi: \Gamma \longrightarrow  H$ 
which satisfies:  
i)
$\psi(\Gamma) \leq  H$ is a discrete, Zariski-dense subgroup, 
ii) $\bH$ has a strong unipotent radical,  
and iii) $\dim U =\rank \Gamma$. 
Let $H'= \bH'_\bbR$ be another real linear algebraic group,
$\psi': \Gamma \longrightarrow H'$ an
embedding of $\Gamma$ into $H'$ so that i) to iii) 
are satisfied with respect to $\bH'$. 
Hence, as a consequence of the rigidity 
property (*), there exists an 
algebraic isomorphism  $\Phi: H \rightarrow H'$ 
so that $\psi' = \Phi \circ \psi$.
We call the solvable 
by finite real linear algebraic group $H_\Gamma:=  H$ 
the real algebraic hull for $\Gamma$. 

\paragraph{The real algebraic hull for $G$}
Let $G$ be a connected, simply connected solvable Lie-group. By 
\cite[Proposition 4.40]{Raghunathan}, 
there exists an algebraic
hull for $G$. This means that 
there exists an $\bbR$-defined linear algebraic
group $\bH_G$, and an injective Lie-homomorphism 
$\psi: G \longrightarrow (\bH_G)_{\bbR}$ 
so that i)' $\psi(G) \leq \bH_G$ is a Zariski-dense subgroup, 
ii)' $\bH_G$ has a strong unipotent radical $\bU$, and iii)' 
$\dim \bU = \dim G$. Moreover, $\bH_G$ satisfies rigidity 
properties analogous to the rigidity properties of $\bH_\Gamma$. 
Let $H_G= (\bH_G)_\bbR$. Then there exists
a continuous injective homomorphism $\psi: G \longrightarrow H_G$
which has Zariski-dense image in $H_G$. 
As for real hulls of discrete groups,
these data are uniquely defined up to composition with an 
isomorphism of real algebraic groups, and
we call $H_G$ the real algebraic hull for $G$.  
We consider henceforth a 
fixed continuous Zariski-dense inclusion $G \leq H_G$.

Let $N$ denote the nilpotent radical of $G$, 
i.e., the maximal, connected
nilpotent normal subgroup of $G$,
and let $U_G$ denote the unipotent
radical of $H_G$. Then (compare 
the proof of Lemma \ref{Fittu}) 
$N \leq U_G=\ur(H_G)$, so that 
$N$ is the connected component of 
$\ur(G)=  G \cap  \ur(\bH_ G)$.
We remark further that $G$ is a normal subgroup of  $H_G$. In fact, 
$N \leq U_G$ is Zariski-closed in $H_G$, and $[G,G] \leq N$ implies therefore
that $[H_G, H_G] \leq N$.
Let $T$ be a maximal $d$-subgroup of $H_G$. 
We consider the decomposition $H_G = U_G \cdot T$. 
Since $H_G$ is decomposed as a product of varieties,
the projection map
$\tau_T: H_G \rightarrow U_G, \;  g= ut \mapsto u $,
onto the first factor of the splitting
is an algebraic morphism.

\begin{proposition} \label{Gsplitting} 
Let $G$ be connected, simply connected solvable Lie group, 
and $G \leq H_G$ a continuous, Zariski-dense 
inclusion into its real algebraic hull. 
Then $G$ is a closed normal subgroup of $H_G$.   
Moreover, if $T \leq H_G$ is a maximal $d$-subgroup 
then 
\begin{equation*} 
 H_G = G  T \;, \; G \cap T = \{1 \} \;\;  . 
\end{equation*} 
Let $U_G = \ur(H_G)$ denote the unipotent radical 
of $H_G$. Then the 
algebraic projection map $\tau_T: H_G \rightarrow U_G$ restricts
to a diffeomorphism $\tau: G \rightarrow U_G$.  
\end{proposition}    
\begin{proof} 
Let $C$ be a Cartan subgroup of $G$. Then $C$ is 
nilpotent, and $G= N C$, where $N \leq U_G$ is the nilradical of $G$. 
Let us put $S= C_s = \{ g_s \mid g \in C \}$,  so that 
$C \leq C_u \times S$. Note that $C_u$ is a closed subgroup of $U_G$, 
and $S$ is an abelian subgroup of $H_G$ which is centralized
by $C$. Let $T \leq H_G$ be a maximal $d$-subgroup which
contains $S$. 
Since $H_{G} = \ac{G} \leq N C_u T$, we conclude
that $U_G = N C_u$ and $H_G= G T$.  
It follows that the crossed 
homomorphism $\tau: G  \rightarrow U_G$ 
is surjective, in fact, since $\dim U_G = \dim G$ it 
is a covering map.
Since $U_G$ is simply connected $\tau$ 
must be a diffeomorphism. 
Therefore $T \cap G = \{1\}$. 
From the above remarks, $G$ is a normal
subgroup of $H_G$. Let $\pi_T: H_G \rightarrow T$
denote the projection map onto the second  
factor of the splitting $H_G= U_G \cdot T$.
Then 
$G= \{  g= u \theta(u) \mid u \in U_G \}$, 
where $\theta= \pi_T \tau^{-1}: U \rightarrow T$ is 
a differentiable map. Therefore $G$ is a closed 
subgroup. 
\end{proof} 

Let $\Gamma \leq G$ be a lattice. We call $\Gamma$ a Zariski-dense
lattice if $\Gamma$ is Zariski-dense in $H_G$. We remark:  

\begin{proposition} \label{algebraicH}
Let $G$ be a connected, simply  connected 
solvable Lie group, and\/ $\Gamma \leq G$ a Zariski-dense lattice. Then
the real algebraic hull $H_{G}$ is a real algebraic hull for $\Gamma$. 
\end{proposition} 
\begin{proof} By the inclusion $\Gamma \leq G$ we have
an inclusion $\Gamma \leq H_G$. 
Since $\Gamma$ is cocompact, $\rank \Gamma = \dim G = \dim \ur(\bH_G)$.
Therefore $\bH_G$ is a $\bbR$-defined algebraic hull for $\Gamma$. 
By the rigidity property (*),
there exists an $\bbR$-defined isomorphism $\bH_\Gamma \rightarrow \bH_G$.
In particular, there is an induced algebraic
isomorphism of the groups of real points 
$H_\Gamma$ and  $H_G$. 
\end{proof} 

Identifiying, $\Aut(G)$ with $\Aut(\lie{g})$, where $\lie{g}$ is
the Lie-algebra of $G$, we obtain a natural structure of real 
linear algebraic group on $\Aut(G)$. Let $\bH$ be a 
solvable by finite linear algebraic group, and let
$\Aut_a(\bH)$ denote its group of algebraic automorphisms.
In \cite{GP1}, 
it is observed that the group $\Aut_a(\bH)$ is itself a  
linear algebraic group if $\bH$ has a strong unipotent
radical. In particular, $\Aut_a(H_G)$, 
the group of algebraic automorphisms of $H_G$,  
inherits a structure of a real linear algebraic group. 
The rigidity 
of the hull $H_G$ induces an extension homomorphism 
$$ {\cal{E}}: \Aut(G) \hookrightarrow \Aut_a(H_G) \; , 
\; \psi \mapsto \Psi \; .  $$

\begin{proposition} \label{AutGisasubgroup}
The extension homomorphism $  {\cal{E}}: 
\Aut(G) \hookrightarrow \Aut_a(H_G)$
identifies the real linear algebraic group $\Aut(G)$ with
a Zariski-closed subgroup of  $\Aut_a(H_G)$
\end{proposition}
\begin{proof} Let $\lie{h}_G$ denote the Lie-Algebra
of $H_G$. From the inclusion $G \leq H_G$, we have
that $\lie{g} \subseteq \lie{h}_G$.
Since $H_G= H_G^0$, it follows from the discussion in \cite{GP1}[\S3]
that the Lie-functor identifies the group $\Aut_a(H_G)$ with a 
Zariski-closed subgroup  $\Aut_a(\lie{h}_G)$ of $\Aut(\lie{h}_G)$. Consider
$$ \Aut_a(\lie{h}_G, \lie{g}) = 
\Aut_a(\lie{h}_G) \cap \{ \varphi  \mid \varphi(\lie{g})\subseteq
\lie{g} \} \; . $$ 
The rigidity property of the hull implies that 
the restriction map $$\Aut_a(\lie{h}_G, \lie{g}) \, \rightarrow \, \Aut(\lie{g})$$ 
is surjective. Since $G \leq H_G$ is Zariski-dense, the restriction map 
is injective as well. This implies that the restriction map 
induces an isomorphism of real linear algebraic groups. Since
the image of $\Aut(G)$ in $\Aut_a(H_G)$ correponds to the
Zariski-closed subgroup $\Aut_a(\lie{h}_G, \lie{g}) \leq \Aut_a(\lie{h}_G)$,
the proposition follows.  
\end{proof}

\subsection{Affine actions by rational maps} \label{ARactions}
Let $G$ be a group.  We view the affine group $\Aff(G)$
as a group of transformations acting on  $G$ by declaring 
$$ \;\; \; \; (g, \phi) \cdot g' = g \, \phi(g')  \;  , \; \;  \;  
\text{where } (g, \phi) \in \Aff(G), \; g' \in G \; . $$
Let $H$ be a solvable by finite real linear 
algebraic group with a strong unipotent 
radical. Let $\Aut_a(H) \leq \Aut(H)$ denote 
its group of algebraic automorphisms. 
We remark that, since $H$ has a strong unipotent
radical,  {\em $\Aut_a(H)$ is a real
linear algebraic group\/} (as follows from  \cite[\S 4]{GP1}),
and so is $$ \Aff_a(H)= H \rtimes \Aut_a(H) \;.  $$
Let $T$ be a maximal $d$-subgroup 
of $H$, and  $U = \ur(H)$. For $h \in H$, let $c(h): H \rightarrow H$
denote the inner automorphism $l  \mapsto hlh^{-1}$ of $H$.
If $L \leq H$ is a normal subgroup $c_L(h)$ denotes the
restriction of $c(h)$ on $L$.   
Let $h= u   t$ be a decomposition of $h \in H$ with 
respect to the algebraic splitting $H = U \cdot T$. 
Then we have  a homomorphism of real algebraic groups
$$ \alpha_T: H \longrightarrow \Aff_a(U),  \; \;  h=ut \mapsto (u , c_U(t))\; \;  .$$ 
Since $H$ has a strong unipotent radical, 
the homomorphism $\alpha_T$ is injective. Similarly, 
if $G \leq H$ is a normal subgroup of $H$, $H= G T$ and
$G \cap T = \{1\}$, we define  
$$ \beta_T:  H \longrightarrow \Aff(G), \; \;  
h=gt \mapsto (g , c_G(t)) \; \; . $$ 

\begin{lemma} \label{equivariant} 
Let $H$ be a solvable by finite 
real algebraic group with a strong unipotent radical $U$, 
and let $T \leq H$ be a maximal $d$-subgroup. Assume 
there exists a connected Lie subgroup $G$ of $H$ which
is normal in $H$, so that
$H= G T$, $H \cap T= \{1\}$. Then $\beta_T: H \rightarrow \Aff(G)$
is an injective continuous homomorphism.     
Moreover, the projection $\tau: G \rightarrow U$, induced 
by the splitting $H= U \cdot T$,  is a diffeomorphism 
which is equivariant with respect to the affine actions  $\beta_T$
and $\alpha_T$.       
\end{lemma} 
\begin{proof}  We remark first that $H= G T$ implies 
that the Zariski-closure $\ac{G} \leq H$ 
contains the unipotent radical $U$ of $H$. 
The argument given in the proof of 
Proposition \ref{Gsplitting} shows that $\tau$ is
a diffeomorphism. Using the notation in the proof 
of Proposition \ref{Gsplitting}, we can write
$$\beta_T(h)= (\, \tau_T(h) \theta(\tau_T(h)) \, , \,  
c_G\!\left( \theta(\tau_T(h))^{-1} \pi_T(h)\right) \, ) \; .$$ 
This shows that $\beta_T$
is continuous. It is also injective: Assume that $\beta_T(gt)=1$.
Then, in particular, $c_G(t) = id_G$. 
Hence,  $t$ centralizes the unipotent radical $U$. 
Since $t$ is semisimple and $H$ has a strong unipotent
radical, this implies $t=1$. Therefore, $h \in G$. But $\beta_T$ is 
clearly injective on $G$, proving that $\beta_T$ is injective.    
Finally, let $h \in H$, $g \in G$. 
Then an elementary calculation shows 
that $\tau(\beta_T(h) \cdot g) = \alpha_T(h) \cdot \tau(g)$,
proving that $\tau$ is equivariant. 
\end{proof} 

Finally, we briefly remark how the affine action $\alpha_T$ 
depends on the choice of maximal $d$-subgroup in $H_\Gamma$. 
Let $T' \leq H_\Gamma$ be another maximal $d$-subgroup.
By Proposition \ref{Usplitting}, there exists $v \in U$ so
that  $T' = v T v^{-1}$. Let $h= ut$, where $u \in U$, $t \in T$. 
The decomposition $h=u't'$ of $h$ relative to $T'$ is 
given by $u'= uv^tv^{-1}$, $t'= t^v$. Hence,  
\begin{lemma} \label{conjugacy}  
Let $R_v: U \rightarrow U$ denote
right-multiplication with $v$ on $U$, and $T' = v T v^{-1}$.
Then, for all $h \in H$, 
$\alpha_T(h) \circ R_v = R_v \circ \alpha_{T'}(h)$. 
\end{lemma} 

\subsection{Syndetic hulls}
The notion of {\em syndetic hull\/} of a solvable subgroup
of a linear group is due to Fried and Goldman, 
cf.\  \cite[\S 1.6 ]{FriedGoldman}.
Fried and Goldman introduced this notion in the
context of affine crystallographic groups. We will employ
the syndetic hull to prove that standard $\Gamma$-manifolds
are infra-solvmanifolds. 
We use the slightly modified 
definition for the syndetic hull which is 
given in \cite{GS}. Let $V$ be a finite-dimensional real 
vector space. 

\begin{definition} \label{syndetic} 
Let $\Gamma$ be a polycyclic subgroup of $\GL(V)$, and
$G$ a closed, connected subgroup of $\GL(V)$ such that $\Gamma \leq G$. 
$G$ is called a
{\em syndetic hull\/} of $\Gamma$ if $\Gamma$ is a Zariski-dense 
(i.e., $G \leq \ac{\Gamma}$) uniform lattice in $G$, 
and $\dim G = \rank \Gamma$. 
\end{definition} 

The syndetic hull for $\Gamma$ is necessarily 
a connected, simply connected solvable Lie group.
If $\Gamma \leq \GL(V)$ is discrete, $\Gamma \leq (\ac{\Gamma})_0$,
and $\Gamma/ \ur(\Gamma)$ is torsion-free then it  
is proved in  \cite[Proposition 4.1, Lemma 4.2]{GS})
that $\Gamma$ has a syndetic hull. In particular, any
virtually  polycyclic linear group has a normal finite
index subgroup which possesses a syndetic hull. 
We need the following slightly refined result:
 
\begin{proposition} \label{syndhull}
Let $H \leq \GL(V)$ be a Zariski-closed subgroup.
Let $\Delta \leq H$ be a virtually polycyclic discrete
subgroup, Zariski-dense in $H$.  
Then there exists a finite index normal subgroup 
$\Gamma_0 \leq \Delta$, and a syndetic hull $G$
for $\Gamma_0$, $\Gamma_0 \leq G \leq H$, so that $G$
is normalized by $\Delta$. 
\end{proposition}   

A similar result is also stated in \cite[\S 1.6 ]{FriedGoldman}.
However, the proof given in \cite{FriedGoldman} is faulty. 
We refine the proof of \cite[Proposition 4.1]{GS} 
a little to obtain Proposition \ref{syndhull}. 
Also we warn the reader that a syndetic hull $\Gamma \leq G$
is (in general) not uniquely determined by $\Gamma$, neither a
{\em good syndetic hull\/} (cf.\ \cite{GS}) is uniquely determined by
$\Gamma$. (See \cite[\S 9]{GS}.) 

\medskip 
\begin{prf}{Proof of Proposition \ref{syndhull}:}
There exists a normal polycyclic subgroup 
$\Gamma \leq \Delta$ of finite index with the
following properties:  
$\Gamma \leq H_{0}$, $[\Gamma,\Gamma] \leq \ur(H)$, and 
$\Gamma \rmod \ur(\Gamma)$ is torsion-free. 
We consider the abelian by finite Lie group $T=H/N$, 
where $N= \ac{[\Gamma,\Gamma]}$.
Let $p: H \rightarrow T$ denote the 
projection, and $\pi: \hat{E} \rightarrow T$ the universal cover. Then 
$\hat{E}$ is an extension of a vector space $E$ by some 
finite group $\mu$. Let $\hat{S}= \pi^{-1}( \Delta N/N)$, 
$S = \hat{S} \cap E$. Then $\mu= \hat{S} /S$. 
Now $K = \ker \pi \subset S$ is invariant
by the induced action of $\mu$ on $S$. From Maschke's theorem 
we deduce that there exists a $\mu$-invariant complement 
$P \subset S$ of $K$ so that $K P$ is of finite index in $S$. 

Now define $\bar{P}$ to be the real vector space 
spanned by $P$, $G = p^{-1}(\pi(\bar{P}))$,  
and $\Gamma_0 = G \cap \Gamma$. Then $\Gamma_0$ is
of finite index in $\Gamma$, and $G$ is a syndetic 
hull for $\Gamma_0$ in the sense of Definition \ref{syndetic}. 
Since $\bar{P}$ is invariant by $\mu$ the Lie group $G$ is 
normalized by $\Delta$. 
\end{prf}

\section{Standard $\Gamma$-manifolds} \label{sstandardG}
Let $\Gamma$ be a torsion-free virtually polycyclic 
group. The purpose of this section is to explain  
the construction of standard $\Gamma$-manifolds 
and to prove Theorem \ref{main1} and Theorem \ref{main}. 

\subsection{Construction of standard $\Gamma$-manifolds}  
Let $H_\Gamma$ be a real algebraic hull for $\Gamma$, 
and fix a Zariski-dense embedding  $\Gamma \leq H_\Gamma$. 
Let $T$ be a maximal $d$-subgroup of $H_\Gamma$, and
put $U= \ur(H_\Gamma)$ for the unipotent radical of $H_\Gamma$.  
We consider the affine action $\alpha_T: H_\Gamma \rightarrow
\Aff_a(U)$ which is defined by the splitting $H_\Gamma = U \cdot T$.
Since $U$ is strong in $H_\Gamma$, the 
homomorphism $\alpha_T$ is injective. Let 
$$ M_{\Gamma,\alpha_T}  \; = \; 
{_{\displaystyle \alpha_T(\Gamma)}}  \backslash \, U $$
denote the quotient space of the affine action
of $\Gamma$ on $U$. We will show that $M_{\Gamma,\alpha_T}$
is a compact manifold with 
fundamental group isomorphic to $\Gamma$. 
In fact, the proof implies that $M_{\Gamma,\alpha_T}$ 
is an infra-solvmanifold. 
We also show
that the diffeomorphism class of $ M_{\Gamma,\alpha_T}$
depends only on $\Gamma$, 
not on the choice of maximal $d$-subgroup
$T$ in $H_\Gamma$, nor on the particular embedding
of $\Gamma$ into $H_\Gamma$. In fact, we show that 
the corresponding actions of $\Gamma$ are affinely
conjugate. We call $M_\Gamma = M_{\Gamma,\alpha_T}$  
a {\em standard $\Gamma$-manifold}.
\medskip 

\begin{prf}{Proof of Theorem \ref{main1}:}
We show first that $M_\Gamma$ is an 
infra-solvmanifold. Let $\Gamma_0$ be a finite
index normal subgroup of $\Gamma$ so that there
exists a syndectic hull $\Gamma_0 \leq G \leq H_\Gamma$ for $\Gamma_0$.
By Proposition \ref{syndhull}, we 
may also assume that $G$ is normalized by $\Gamma$.
By the defining properties of the hull, $(H_\Gamma)^0$ is 
a real algebraic hull for $G$.
Since $G$ is a normal subgroup in its hull, it follows that
$G$ is a normal subgroup of $H_\Gamma= \Gamma(H_\Gamma)^0$. 
Let $T$ be a maximal $d$-subgroup of $H_\Gamma$. 
We infer from 
Proposition \ref{Gsplitting} that  
$H_\Gamma = G T$, $G \cap T = \{1\}$. 
Let $\beta_T: H_\Gamma \rightarrow \Aff(G)$
denote the affine action which is defined 
by this splitting.
Lemma \ref{equivariant} implies that
the affine action $\beta_T$ is effective. 
Note that $\beta_T(\Gamma) \cap G$ contains $\Gamma_0$,
hence $\beta_T(\Gamma)$ is discrete in $\Aff(G)$ and 
$hol(\beta_T(\Gamma)) \leq \Aut(G)$ is finite. 
Therefore, the quotient space  
$$ M_{\beta_T} \; = \; {_{\displaystyle \beta_T(\Gamma)}} \backslash\, G$$
is an infra-solvmanifold. Since $G$ is diffeomorphic to $\bbR^n$, 
the fundamental group $\pi_1( M_{\beta_T})$
is isomorphic to $\Gamma$. Since $G$ is a syndetic hull,
$M_{\beta_T}$ is compact.   
Let  $\tau: G \rightarrow U$
be the projection map which is induced 
by the splitting $H_\Gamma= U \cdot T$. 
By Lemma \ref{equivariant}, $\tau$ induces 
a diffeomorphism $\bar{\tau}: \beta_T(\Gamma) \backslash G
\rightarrow \alpha_T(\Gamma) \backslash U$. 
Hence 
$\alpha_T(\Gamma) \backslash U$ is 
diffeomorphic to a compact infra-solvmanifold. 

Note that the diffeomorphism class of  $M_\Gamma$
does not depend on the choice of maximal $d$-subgroup in $H_\Gamma$. 
In fact, let $T' \leq H_\Gamma$ be another maximal $d$-subgroup.
Then, by Proposition \ref{Usplitting}, there exists $v \in U$ so
that $T' = v T v^{-1}$. By Lemma \ref{conjugacy}, 
$R_v: U \rightarrow U$ induces a smooth diffeomorphism 
$M_{\Gamma,\alpha_T} \rightarrow M_{\Gamma,\alpha_{T'}}$. 
The diffeomorphism class of $M_\Gamma$ is also independent of
the particular choice of Zariski-dense embedding of $\Gamma$ into $H_\Gamma$. 
Let $\Gamma' \leq H_\Gamma$ be a Zariski dense
subgroup isomorphic to $\Gamma$, and let
$\phi: \Gamma \rightarrow \Gamma'$  
be an isomorphism.
By the rigidity of the real algebraic hull, 
there exists an algebraic
automorphism 
$\bar{\phi}: H_\Gamma \rightarrow H_{\Gamma}$
extending $\phi$. 
The restriction of $\bar{\phi}$ on the unipotent
radical $U$ of $\bH_\Gamma$ projects 
to a diffeomorphism $M_{\Gamma,\alpha_T} \rightarrow 
M_{\Gamma',\alpha_{\bar{\phi}(T)}}$ 
which induces $\phi$ on the level of fundamental groups.
Similarly, any given automorphism of $\pi_1(M_{\Gamma,\alpha_T})$
corresponds to an automorphism
$\phi$ of $\Gamma$. The algebraic extension 
$\bar{\phi}$ projects to  a
diffeomorphism of $M_{\Gamma,\alpha_T}$ inducing $\phi$
on $\pi_1(M_{\Gamma,\alpha_T})$. 
\end{prf}  

We also remark: 
\begin{proposition} \label{storusrank}
Every standard $\Gamma$-manifold
$M_\Gamma$ admits a smooth effective action of an 
$r$-dimensional torus $T^r$, where $r = \rank \Z(\Gamma)$.
\end{proposition}
\begin{proof} Let $U_Z = \ac{\Z(\Gamma)}$ be the Zariski-closure 
of $\Z(\Gamma)$ in $H_\Gamma$. Since $U_Z \leq \Z(H_\Gamma)$, 
$U_Z \leq U$, and $\dim U_Z= r$. It follows that $\alpha_T(U_Z)$ induces
a free maximal torus action on $M_\Gamma$. 
\end{proof}
Theorem \ref{main2} and Proposition \ref{storusrank}
imply Corollary \ref{torusrank}.

\subsection{Affine actions on unipotent groups}
Here we show that every compact manifold 
which arises by (solvable by finite) affine actions on unipotent
groups is diffeomorphic to a standard $\Gamma$-manifold. 

\begin{proposition}
Let $U$ be a connected, simply connected nilpotent
Lie group. Let $\Delta \leq \Aff(U)$ be a solvable by finite 
subgroup which acts freely and properly on $U$ with compact
quotient manifold $M = \Delta \backslash U$. 
Let $T \leq \ac{\Delta} \leq \Aff(U)$ be a 
maximal $d$-subgroup. Then  $M$ is diffeomorphic to 
$\alpha_T(\Delta) \backslash \ur(\ac{\Delta})$. 
\end{proposition}
\begin{proof} 
We decompose $\ac{\Delta} = \ur(\ac{\Delta}) T$. 
Since any maximal $d$-subgroup of
$\Aut(U)$ is maximal in $\Aff(U)$, we may assume
(after conjugation of $\Delta$ with a suitable element of 
$\Aff(U)$) that $T \leq \Aut(U)$. In particular, 
$\ac{\Delta} \cdot 1= \ur(\ac{\Delta}) \cdot 1$, 
and hence $\Delta$ acts on the orbit $O= \ur(\ac{\Delta}) \cdot 1 \leq U$. 
Since $O$ is the homogeneous space of a connected unipotent 
group acting on $U$, it is a submanifold diffeomorphic 
to $\bbR^k$, $k \leq \dim U$. The connected, simply connected 
solvable Lie group $\Delta_0$ acts freely on $O$, 
and the quotient space $\Delta_0 \backslash O$ is
a simply connected aspherical manifold. Hence, the
quotient space $M_O = \Delta \backslash O$ is a
manifold with fundamental group isomorphic 
to $\Gamma$, and homotopy equivalent to an Eilenberg-Mac Lane 
space $K(\Gamma,1)$. Since $M$ is an aspherical 
compact manifold with fundamental group $\Gamma$, its
dimension equals the cohomological dimension of $\Gamma$. 
This implies that $\dim M \leq \dim M_O$, and 
consequently $O= U$. In particular, $\ur(\ac{\Delta})$
acts transitively on $U$ and the orbit map
$$ o: \ur(\ac{\Delta}) \rightarrow U \; , \; \delta \mapsto \delta \cdot 1$$  
in $1 \in U$ is a diffeomorphism. 
Using $T \leq \Aut(U)$, 
it is straightforward to verify 
that $o$ is $\Delta$-equivariant
with respect to the affine action $\alpha_T$ of
$\Delta$ on $\ur(\ac{\Delta})$.  
Hence, $M$ is diffeomorphic to 
$\alpha_T(\Delta) \backslash \ur(\ac{\Delta})$. 
\end{proof}

The next result shows how the algebraic 
hull enters the picture:

\begin{proposition} \label{Udiscrete}
Let $\Gamma \leq \Aff(U)$ be virtually polycyclic,
such that $\Gamma$ acts freely and properly discontinuously
on $U$, and with compact quotient $M = \Gamma \backslash U$. 
Then the Zariski-closure $\ac{\Gamma} \leq \Aff(U)$ is an 
algebraic hull for $\Gamma$. In particular, $M$
is diffeomorphic to a standard $\Gamma$-manifold. 
\end{proposition}
\begin{proof} Put $H=\ac{\Gamma} \leq \Aff(U)$. 
By the previous proposition, the orbit map 
$o: \ur(H) \rightarrow U$ in $1 \in U$
is a diffeomorphism. In particular, 
$\dim \ur(H) = \rank \Gamma$. Let $T \leq H$
be a maximal $d$-subgroup. We may assume that
$T \leq \Aut(U)$. This shows that 
$T \cap \Z_H(U) = \{1\}$. 
Hence, $H$ has a strong unipotent radical.
It follows that $H$ is an algebraic hull 
for $\Gamma$. Since $o$ is equivariant with
respect to the action $\alpha_T$, $M$
is diffeomorphic to a standard $\Gamma$-manifold.
\end{proof}

Next we consider affine actions which arise from splittings of solvable
by finite linear algebraic groups.  

\begin{proposition} Let $H= U T$ be a solvable linear 
algebraic group.  Let $\Theta \leq H$ be a Zariski-dense 
subgroup such that $\Delta = \alpha_T(\Theta) \leq \Aff(U)$ 
acts freely and properly 
on $U$ with compact quotient manifold $M = \Delta \backslash U$. 
Then $\Gamma = \Delta / \Delta_0$
is virtually polycyclic and $M$ is diffeomorphic to a standard
$\Gamma$-manifold.   
\end{proposition}
\begin{proof} Put $U_{\Delta_0} = \ur(\ac{\Delta_0})$, and
remark that $U_{\Delta_0} \leq U$
under the natural inclusion $U \leq \Aff(U)$.
Since $\Delta_0$ is normal in $\Delta$,
$U_{\Delta_0}$ is normal in $\ac{\Delta}$.
From  $\ac{\Delta_0} \cdot 1= U_{\Delta_0} \cdot 1$
and $\dim \ur(\ac{\Delta_0}) \leq \dim \Delta_0$, we
deduce that $ \Delta_0 \cdot 1 = U_{\Delta_0} \cdot 1$
and also that $\dim  \Delta_0 = \dim U_{\Delta_0}$.  
Let $h \in \Theta$, such that $\alpha_T(h) \in \Delta_0$.
Then  $h = u_h t$, where $u_h \in U_{\Delta_0}$ 
and $t \in T$. Moreover, 
$$\alpha_T(h) \cdot u = u_h u^t = u_h (t u t^{-1} u^{-1}) u \; .$$ 
Remark that $\ac{\Delta_0}$, and $U_{\Delta_0}$ 
are normal in $\alpha_T(H)$. This implies that   
$h ^u = v t$, where $v \in U_{\Delta_0}$, 
$v=  (u u_h u^{-1}) (utu^{-1}t^{-1})$.
Furthermore $u u_h u^{-1} \in U_{\Delta_0}$, and consequently 
$ut u^{-1}t^{-1} \in U_{\Delta_0}$.  
Since $\Delta_0$ acts freely, this implies that 
$\Delta_0 \cdot u = \ur(\ac{\Delta_0}) \cdot u$.
Hence, $\Delta_0$ and $U_{\Delta_0}$
have the same  orbits on $U$.

 
Put $L  = H/U_{\Delta_0}$, 
and $U_L=  U /U_{\Delta_0}$. 
Let $\pi: H \rightarrow L$ be the
quotient homomorphism. Put 
$\Upsilon:= \pi(\Theta)$. 
Then $\Upsilon$ 
is a Zariski-dense subgroup of $L$. 
Decompose $L = U_L \pi(T)$. Evidently, $\pi$
induces a diffeomorphism of quotient spaces 
$\alpha_{T}(\Theta) \backslash U 
\rightarrow \alpha_{\pi(T)}(\Upsilon) \backslash U_L$. 
In particular, $M$ is diffeomorphic to 
$\alpha_{\pi(T)}(\Upsilon) \backslash U_L$.
Moreover, $\alpha_{\pi(T)}(\Upsilon)$ is 
a discrete solvable subgroup of $\Aff(U_L)$ and 
isomorphic to $\Delta/ \Delta_0$. Since $\Aff(U_L)$ has only 
finitely many connected components, a theorem 
of Mostow \cite{Mostow3} implies that  
$\Gamma =\Delta/ \Delta_0$ is virtually polycylic.
By Proposition \ref{Udiscrete}, $M$ is diffeomorphic to a 
standard $\Gamma$-manifold.
\end{proof}
 
Putting the results together, we proved: 

\begin{theorem} \label{Urigidity}
Let $U$ be connected, simply connected nilpotent
Lie group. Let $\Delta \leq \Aff(U)$ be a solvable by finite 
subgroup which acts freely and properly on $U$ with compact
quotient manifold $M = \Delta \backslash U$. 
Then $\Gamma = \Delta/\Delta_0$
is virtually polycyclic, and $M$ is diffeomorphic 
to a standard $\Gamma$-manifold. 
\end{theorem}

\subsection{Rigidity of reductive affine actions}
Let $G$ be a connected, simply connected
solvable Lie group and $\Delta \leq \Aff(G)$
a solvable by finite subgroup which acts
on $G$. Let $H_G$ be an algebraic hull 
for $G$, and fix a Zariski-dense continuous inclusion 
$G \leq H_G$. 
By the rigidity of the hull, there are 
induced inclusions $hol(\Delta) \leq \Aut_a(H_G)$, 
and $\Delta \leq \Aff_a(H_G)$.  
Let $T \leq H_G$ be a maximal $d$-subgroup, and
$U_G$ the unipotent radical of $H_G$. 
Then, $G$ acts
affinely on $U_G$ via the action $\alpha_T$,
c.f.\ \S  \ref{ARactions}.
Note that the orbit map  of this action in 
$1 \in U_G$, $o_T: G \rightarrow U_G$, coincides with
the projection diffeomorphism $\tau: G \rightarrow U_G$. 
Via $o_T$, the affine action of $\Delta$ on $G$
induces then a diffeomorphic action of $\Delta$ 
on $ U_G$. 

\begin{lemma} \label{DactsonU}
Suppose that $hol(\Delta) \leq \Aut_a(H_G)$ stabilizes $T$.
Then the action of  $\Delta$ on $ U_G$ induced by the
orbit map $o_T: G \rightarrow U_G$ is affine. 
\end{lemma}
\begin{proof} A straightforward computation shows 
that the lemma is true:  Let $\delta = (h, \phi) \in \Aff(G)$,
where $h \in G$, $\phi \in hol(\Delta)$. We consider $\phi$ 
henceforth as an element of $\Aut(H_G)$. Write $h= u_h t_h$, where
$u_h \in U_G$, $t_h \in T$. Analogously, 
write $g= u_g t_g$, for $g \in G$.
Now, $\delta \cdot g  = h \phi(g) = u_h t_h \phi(u_g) \phi(t_g)$.
By our assumption, $\phi(t_g) \in T$. Hence, 
\begin{equation*}\begin{split}
 o_T(\delta \cdot g) &=  \tau( u_h \phi(u_g)^{t_h}  t_h \phi(t_g)) = 
 u_h  \phi(u_g)^{t_h} \\
 &= \alpha_T(h) \cdot \phi(u_g) = (\alpha_T(h) \circ \phi) \cdot o_T(g) \; . 
\end{split} \end{equation*}
Therefore, the action of $\delta$ on $G$, corresponds 
to the action of $\alpha_T(h) \circ \phi$ on $U_G$.  
\end{proof}

\begin{lemma} \label{LT}
Let $L \leq \Aut_a(H_G)$ be a reductive subgroup. 
Then $L$ stabilizes a maximal torus $T \leq H_G$.
\end{lemma}
\begin{proof} Consider the semi-direct product 
$H_L = H_G \rtimes L$. Then 
$\ur(H_L) = U_G$ is the unipotent
radical of $H_L$. Let $S$ be a maximal reductive subgroup in
$H_L$ which contains $L$ such that $H_L = U_G S$. 
Then $T = S \cap H_G$ is a $d$-subgroup in $H_G$
which is normalized by $L$, and 
$S= T L$. The latter equality  shows  
that $H_G= U_G T$ and therefore 
$T$ is a maximal torus of $H_G$. 
\end{proof}

In the light of Theorem \ref{Urigidity}, Lemma \ref{DactsonU} and 
Lemma \ref{LT} prove the following:  

\begin{theorem} \label{Grigidity} 
Let $\Delta \leq \Aff(G)$ act freely on $M$
with quotient space $\Delta \backslash G$ a compact manifold.  
Assume further that  
$\hol(\Delta) \leq \Aut_a(H_G)$ 
is contained in a reductive subgroup of $\Aut_a(H_G)$.
Then $M= \Delta \backslash G$ is diffeomorphic to 
a standard $\Gamma$-manifold. 
\end{theorem}

Theorem \ref{Grigidity} implies part ii) of
Theorem \ref{main}. In fact, by Proposition 
\ref{AutGisasubgroup} the assumption ii) implies 
that $\hol(\Delta)$ is contained in a reductive
subgroup of $\Aut_a(H_G)$. In particular, 
condition ii) is satisfied if the Zariski-closure of
$\hol(\Delta)$ in $\Aut(G)$ is compact. This
proves then Theorem \ref{main2}.

\section{Geometry of infra-solvmanifolds}
We derive a few consequences of our 
proof which concern the existence and uniqueness 
of certain geometric structures on infra-\-solv\-mani\-folds.   
As another application we construct a finite-dimensional 
complex which computes the cohomology
of an infra-\-solv\-mani\-fold. 

\subsection{Infra-solv geometry}
Let $M$ be a compact infra-solvmanifold.
A pair $(G, \Delta)$, $\Delta \subset \Aff(G)$,
so that $M$ is diffeomorphic to $\Delta \backslash G$ is called 
a presentation for $M$. By the proof of
Theorem \ref{main1}, every standard $\Gamma$-manifold 
admits a presentation 
$(G, \Gamma)$ so that $\Gamma \leq \Aff(G)$ 
is discrete with finite holonomy group $\hol(\Gamma)$. 
Hence, by Theorem \ref{main2}, every 
compact infra-solvmanifold has such a
presentation. (The appendix of \cite{FJ2} is devoted 
to proving that every infra-solvmanifold has a
presentation with finite holonomy.) 

\begin{corollary} Every compact infra-solvmanifold
$M$ admits a discrete presentation with finite holonomy.  
\end{corollary}

Let $(G, \Gamma)$ be a discrete presentation for $M$
with finite holonomy.  Then  
$M$ is finitely covered by the 
homogeneous space $\Gamma \cap G \backslash G$ 
of the solvable Lie-group $G$. The group $G$ carries a natural flat 
(but not necessarily torsion-free) 
left invariant connection which is
preserved by $\Aff(G)$. Since the 
presentation is discrete, $M$
has a flat connection inherited
from $G$. The group of
covering transformations of $\Gamma \cap G \backslash G \rightarrow M$
is acting by connection preserving 
diffeomorphisms. This geometric property 
distinguishes infra-solvmanifolds
from the larger class of compact manifolds 
which admit a finite covering by a 
solv-manifold. (Compare also \cite{Tuschmann}
for a similar discussion.)
One should note however that the 
Lie group $G$, and discrete
presentation $(G,\Gamma)$ is not uniquely determined by 
$M$. However, Wilking
\cite{Wilking} proved that 
every infra-solv\-ma\-ni\-fold
is modelled in a canonical way 
on an affine isometric action on 
a super-solvable Lie-group.

Our approach implies that, 
dropping the condition of finite holonomy, 
there is a canonical choice of flat geometry on $M$ 
which is modelled on a nilpotent Lie group. Let $U$ be
a unipotent real algebraic group, $\Gamma \leq
\Aff(U)$ a discrete subgroup which acts properly
discontinuously on $U$. (It is not 
required that the holonomy of $\Gamma$ be 
finite.) Then $\Gamma$ preserves
the natural flat invariant connection on 
$U$, and there is an induced
flat connection on the quotient manifold $M$.
We say that $M$ has an {\em affinely flat geometry
modelled on $U$}. Let $U_\Gamma$
denote the unipotent radical of the real algebraic hull   
of $\Gamma$. We call $U_\Gamma$ the {\em unipotent 
shadow of $\Gamma$}. 

\begin{corollary} Every compact 
infra-solvmanifold $M$ admits an affinely
flat geometry modelled on the simply connected 
nilpotent Lie group $U_{\pi_1(M)}$.  
\end{corollary}

\paragraph{Toral affine actions}
A natural question is the following: 
\emph{
Given $G$ a simply connected solvable Lie group. 
Which polycyclic groups act affinely on $G$ with a
compact quotient manifold?}
In the particular
case where $G$ is abelian, 
this question asks for the  classification of affine 
crystallographic groups. This is a well known and
difficult geometric problem. (Compare \cite{Baues1}, 
and also the references cited therein for some 
recent results). 

Some
answers to the above question 
can be given when putting restrictions
on the holonomy. 
Let us call the holonomy $hol(\Gamma) \leq \Aut(G)$ 
\emph{toral} if the Zariski-closure of $hol(\Gamma)$ 
is a reductive subgroup of $\Aut(G)$. In particular,
$hol(\Gamma)$ is toral if its closure is 
compact. Hence, infra-solvmanifolds come from
toral actions. Note that also standard $\Gamma$-manifolds are
constructed using toral affine actions. 
Now let $U_G$ denote the unipotent
radical of $H_G$, and $U_\Gamma$ the radical
of $H_\Gamma$. We remark:

\begin{proposition} Let $\Gamma$ be torsion-free
virtually polycylic, acting on the connected, 
simply connected solvable Lie group $G$
with compact quotient space and toral holonomy.
Then $U_G = U_\Gamma$. 
\end{proposition}
\begin{proof} Let $T_1$ denote the
Zariski-closure of $\hol(\Delta)$ in $\Aut_a(H_G)$. 
By the assumption, $T_1 \leq \Aut_a(H_G)$ is a $d$-subgroup. 
By Proposition \ref{AutGisasubgroup}, $T_1$ stabilizes 
$G$, i.e., $T_1 \leq \Aut(G)$. 
Also $T_1$ stabilizes a maximal torus $T \leq H_G$. 
Then the corresponding projection map $\tau: G \rightarrow U_G$ induces
an affine action of $\Gamma$ on $U_G$. By the proof 
of Lemma \ref{DactsonU}, the image of $\Gamma$ in $\Aff(U_G)$ is
contained in $U_G \rtimes T \, T_1 \leq \Aff(U_G)$. 
By Proposition \ref{Udiscrete},  the Zariski-closure 
of $\Gamma$ in $U_G \rtimes T\, T_1$ is an 
algebraic hull for $\Gamma$. 
Hence, $U_G =  U_\Gamma$. 
\end{proof}
In particular, if $\Gamma$ acts isometrically
on $G$, then $U_G = U_\Gamma$.  

\subsection{Polynomial geometry}
The construction of standard $\Gamma$-manifolds 
was carried out in the category of real algebraic groups. 
In fact, as noted above such a manifold is obtained as a quotient
space $\Gamma \backslash U$, where $U$ is a unipotent
real algebraic group, and $\Gamma \leq \Aff_a(U)$ 
is a properly discontinuous subgroup. 
In particular, $\Gamma$ acts algebraically
on $U$. A differentiable map of $\bbR^n$
is called a polynomial map if its coordinate
functions are polynomials. A polynomial diffeomorphism is 
a polynomial map which has a polynomial inverse.
A group of polynomial 
diffeomorphisms of $\bbR^n$ is called bounded
if there is a common bound for the degrees of 
the polynomials which describe its elements. 
It is known that any algebraic group action
on $\bbR^n$ is bounded. 
Now, since $U$ is diffeomorphic 
to $\bbR^n$, $n= \dim U$, and since the diffeomorphism 
is given by the exponential map $\exp: {\lie u}= \bbR^n \rightarrow U$
which actually is an algebraic map, we obtain:

\begin{corollary} Every torsion-free 
virtually polycyclic group $\Gamma$
acts faithfully as a discrete group of 
bounded polynomial diffeomorphisms
on $\bbR^n$, $n = \rank \Gamma$. The
quotient space $\Gamma \backslash \bbR^n$
is diffeomorphic to a standard $\Gamma$-manifold.  
\end{corollary}

Slightly more general, our proof works for all
virtually polycyclic groups which do not have 
finite normal subgroups (see the Appendix). The
existence of such actions was shown previously in \cite{IgDek} by 
different methods. Recently, it was proved 
\cite{BenDek} that a bounded polynomial action of $\Gamma$ on 
$\bbR^n$ is unique up to conjugation by a bounded 
polynomial diffeomorphism. 
Therefore Theorem \ref{main2} implies also the following 
characterization of infra-solvmanifolds. 

\begin{corollary} 
Let $M$ be a compact differentiable 
manifold, aspherical and with a virtually polycyclic 
fundamental group. Then $M$ is diffeomorphic 
to an infra-solv manifold if and only if $M$ 
is diffeomorphic to a quotient space
of\/ $\bbR^n$ by a properly discontinuous bounded group 
of polynomial diffeomorphisms. 
\end{corollary} 

\subsection{Polynomial cohomology}
Let $M$ be an infra-solvmanifold and fix 
a diffeomorphism of $M$ to a quotient space
$\bbR^n/\Gamma$, where $\Gamma$ acts as 
a group of bounded polynomial diffeomorphisms. 
Recall that a differential form (more generally
a tensor field) on $\bbR^n$ is called polynomial
if its component functions relative to the 
standard coordinate system are polynomials. 
Since $\Gamma$ acts by polynomial maps, 
the notion of polynomial differential
form on $M$ is well defined. This gives a subcomplex
$$ \Omega^*_{poly}(M) \, \subset\, \Omega^*(M)$$ 
of the $C^\infty$-de Rham complex $\Omega^*(M)$. The following result 
generalizes a theorem of Goldman \cite{Goldman}
on the cohomology of compact complete affine manifolds:

\begin{theorem} \label{cohomology2}
The induced map on cohomology
$H^*_{poly}(M) \rightarrow H^*(M)$ 
is an isomorphism. 
\end{theorem}
\begin{proof} The idea of the proof which is 
given in \cite{Goldman}
carries over to our situation. 
We pick up the notation of \S \ref{sstandardG}.
Let $\Gamma \leq H_{\Gamma}$, where $H_\Gamma$ is
the algebraic hull of $\Gamma$.  
As explained in section \ref{sstandardG},
$H_{\Gamma}$ acts as as 
a subgroup of $\Aff_a(U)$ on $U$.
Via $\exp: \lu \rightarrow U$,
$H_{\Gamma}$ acts by polynomial maps on $\bbR^n=U$. 
The cohomology 
of $M$ is computed by the complex  $\Omega^*(\bbR^n)^\Gamma$
of $\Gamma$-invariant differential forms on $\bbR^n$. 
Therefore, we have to show that the inclusion of complexes 
$$ \Omega^*_{poly}(\bbR^n)^\Gamma \, \rightarrow \, \Omega^*(\bbR^n)^\Gamma $$
induces an isomorphism on cohomology.
Let $\Gamma_0$ be a finite index normal subgroup 
of $\Gamma$ with syndetic hull $G$, so that 
$$   \Gamma_0 \leq G \leq H= (H_{\Gamma})_0 \; . $$  
We consider now the following inclusion maps of complexes  
$$ \Omega^*(\bbR^n)^{H} 
\rightarrow \Omega^*(\bbR^n)^G  \rightarrow \Omega^*(\bbR^n)^{\Gamma_0} \; \; .
$$
Decompose $H = U S = G S$, where $S$ is a maximal
$d$-subgroup of $H$. By Lemma \ref{equivariant}, 
$G$ acts simply transitively on $\bbR^n$ 
via the affine action $\alpha_T$ of 
$H_{\Gamma}$ on $U$.  
Hence, the complex $\Omega^*(\bbR^n)^G$ identifies 
with the complex $H^*(\lg)$ of left invariant differential forms
on $G$. The action of $G$ on $H^*(\lg)$ which is induced by
conjugation is trivial. Since $G$ is Zariski-dense in $H$, 
$H$ acts trivially on $H^*(\lg)$. 
The affine action
of $S \leq H$ on $U$ corresponds to 
conjugation on $G$.   
It follows that $S$ acts trivially on the 
cohomology of the complex $\Omega^*(\bbR^n)^G$. 
Since $S$ acts reductively on $\Omega^*(\bbR^n)^G$, 
this implies $$ H^*(\Omega^*(\bbR^n)^G) =  H^*(\Omega^*(\bbR^n)^G)^S
= H^*(\Omega^*(\bbR^n)^H) \; .$$
In particular, $\Omega^*(\bbR^n)^{H} 
\rightarrow \Omega^*(\bbR^n)^G$ 
induces an isomorphism on cohomology. 
The map $\Omega^*(\bbR^n)^G  \rightarrow \Omega^*(\bbR^n)^{\Gamma_0}$
is an isomorphism on the cohomology level by a theorem of Mostow
\cite{Mostow5} (see also \cite{Raghunathan}[Corollary 7.29]). 
Hence, the induced map 
$H^*(\Omega^*(\bbR^n)^H) \rightarrow 
H^*(\Omega^*(\bbR^n)^{\Gamma_0})$ is an isomorphism. 
 
Next remark that $\Omega^*(\bbR^n)^{H}= \Omega^*_{poly}(\bbR^n)^H=
\Omega^*_{poly}(\bbR^n)^{\Gamma_0}$.
The first equality follows since every left 
invariant differential form on $U$ is 
polynomial relative to the coordinates given by 
the exponential map.
The second equality follows since $\Gamma_0$ is 
Zariski-dense in $H$. 
We conclude that the natural map $$ H^*(\Omega^*_{poly}(\bbR^n)^{\Gamma_0})
 \rightarrow H^*(\Omega^*(\bbR^n)^{\Gamma_0})$$  is an
isomorphism. This proves that $H^*_{poly}(M_{\Gamma_0}) =
H^*(M_{\Gamma_0})$.

Put $\mu= \Gamma/\Gamma_0$. Then $\mu$ acts on the 
cohomology of $\Omega^*(\bbR^n)^{\Gamma_0}$.
The inclusion map $\Omega^*(\bbR^n)^{\Gamma} \rightarrow \Omega^*(\bbR^n)^{\Gamma_0}$
induces an isomorphism on cohomology 
$$ H^*(M_{\Gamma})\, \rightarrow \, H^*(M_{\Gamma_0})^\mu \; .$$
Similarly, $H^*_{poly}(M_{\Gamma_0})^\mu$ is
isomorphic to the cohomology of the $\mu$-invariant forms 
in $\Omega^*_{poly}(\bbR^n)^{\Gamma_0}$,  implying that
$H^*_{poly}(M_\Gamma) = H^*_{poly}(M_{\Gamma_0})^\mu$.
Hence, $$ H^*_{poly}(M_{\Gamma})= H^*(M_{\Gamma}) \; .$$
The theorem follows.  
\end{proof} 

\begin{prf}{Proof of Theorem \ref{cohomology1}.}
By the previous theorem, $H^*(M_{\Gamma})= H^*_{poly}(M_{\Gamma})$.
Now 
$$ H^*_{poly}(M_{\Gamma}) =  H^*(\Omega^*_{poly}(\bbR^n)^{H_{\Gamma}}) \; .$$
Let $T$ be a maximal $d$-subgroup of $H_\Gamma$.
Since $U$ acts simply transitively on $\bbR^n$, the 
complex $\Omega^*_{poly}(\bbR^n)^{H_{\Gamma}}$ is isomorphic to 
the left-invariant forms on $U$ which 
are fixed by $T$. Let $\lu$ denote the Lie 
algebra of $U$. Since $T$ acts reductively on 
the complex $\Omega^*_{poly}(\bbR^n)^U$, it
follows that 
$$ H^*_{poly}(M_{\Gamma}) =  H^*((\Omega^*_{poly}(\bbR^n)^U)^T)=H^*(\lu)^T \; .$$
\end{prf}


\begin{appendix}


\section{Algebraic hulls for virtually polycyclic groups}
\label{appendixA} 
Let $k \leq \bbC$ be a subfield. 
A group $\bG$ is called a $k$-defined 
linear algebraic group if it is a Zariski-closed subgroup of
$\GL_n(\bbC)$ which is defined by polynomials with 
coefficients in $k$. A morphism of algebraic groups 
is a morphism of algebraic varieties which is also
a group homomorphism. A morphism is defined over $k$
if the polynomials which define it have coefficients
in $k$. It is called a $k$-isomorphism if its inverse
exists and is a morphism defined over $k$. 
Let $\bU = \ur(\bG)$ 
denote the unipotent radical of $\bG$. We say that 
$\bG$ has a {\em strong unipotent 
radical\/} if the centralizer $\C_{\bG}(\bU)$ is contained 
in $\bU$.   

\subsection{The algebraic hull}
Let $\Gamma$ be a virtually polycyclic group. Its maximal nilpotent
normal subgroup $\Fitt(\Gamma)$ is called the Fitting subgroup of
$\Gamma$. Now assume that $\Fitt(\Gamma)$ is 
torsion-free and $\C_\Gamma(\Fitt(\Gamma)) \leq \Fitt(\Gamma)$.
In this case,  
we say that $\Gamma$ has a {\em strong Fitting subgroup}. We 
remark (see also Corollary \ref{strongFitting}) that  
this condition is equivalent
to the requirement that $\Gamma$ has no non-trivial finite normal
subgroups. The following result was announced in \cite{Baues1}.
 
\begin{theorem} \label{alghull} 
Let $\Gamma$ be a virtually polycyclic group with a 
strong Fitting subgroup. 
Then there exists a $\bbQ$-defined linear algebraic group $\bH$ and
an injective group homomorphism $\psi: \Gamma \longrightarrow \bH_{\bbQ}$ 
so that,
\begin{itemize}
\item[i)] $\psi(\Gamma)$ is Zariski-dense in $\bH$, 
\item[ii)] $\bH$ has a strong unipotent radical $\bU$,
\item[iii)] $\dim \bU = \rank \Gamma$. 
\end{itemize} 
Moreover, $\psi(\Gamma) \cap \bH_{\bbZ}$ is of finite index in 
$\psi(\Gamma)$. 
\end{theorem} 

We remark that the group $\bH$ is determined by the conditions
i)-iii) up to $\bbQ$-isomorphism of algebraic groups:
 
\begin{proposition} \label{rigidity1} Let $\bH'$ be a $\bbQ$-defined 
linear algebraic group, $\psi': \Gamma \longrightarrow \bH'_{\bbQ}$ 
an injective homomorphism  which satisfies 
i) to  iii) above. Then there 
exists a $\bbQ$-defined
isomorphism $\Phi: \bH \rightarrow \bH'$ so that 
$\psi' = \Phi \circ \psi$.   
\end{proposition}

\begin{corollary} \label{rigidity2} 
The algebraic hull\/ $\bH_{\Gamma}$ 
of\/ $\Gamma$ is unique up to $\bbQ$-isomorphism
of algebraic groups. 
In particular,  every automorphism $\phi$ of\/ $\Gamma$ 
extends uniquely to a
$\bbQ$-defined automorphism 
$\Phi$ of\/ $\bH_{\Gamma}$. 
\end{corollary}

We call the $\bbQ$-defined linear algebraic group $\bH$ 
the {\em algebraic hull for $\Gamma$}. 
If $\Gamma$ is finitely generated torsion-free nilpotent then
$\bH$ is unipotent and Theorem \ref{alghull} and Proposition
\ref{rigidity1} are essentially
due to Malcev. If $\Gamma$ is torsion-free polycyclic,
Theorem \ref{alghull} is due to Mostow \cite{Mostow2} 
(see also \cite[\S IV, p.74]{Raghunathan} for a different
proof).

\subsection{Construction of the algebraic hull}
Let $\Delta$ be a virtually polycyclic group with $\Fitt(\Delta)$ 
torsion-free. Since $\Delta$ is virtually polycyclic it 
contains a torsion-free polycyclic subgroup $\Gamma$ 
which is normal and of finite index. By Mostow's
theorem there exists an algebraic hull 
$\psi_{\Gamma}: \Gamma \rightarrow \bH_{\Gamma}$
for $\Gamma$. Hence, in particular $\Gamma$ is realized 
as a subgroup of a linear algebraic group.

\paragraph{Embedding of finite extensions} 
We use a standard
induction procedure to realize the finite 
extension group $\Delta$ of $\Gamma$ 
in a linear algebraic group
which finitely extends $\bH_{\Gamma}$.
The procedure is summarized in the next lemma.  

\begin{lemma} \label{induction}
Let $\bG$ be a linear algebraic group
defined over $\bbQ$ and $\Gamma \leq \bG_\bbQ$ a 
subgroup. Let $\Delta$ be a finite extension group of\/
$\Gamma$ so that $\Gamma$ is normal in $\Delta$. 
Let $\Delta= \Gamma r_1  \cup \, \cdots \, \cup \, \Gamma r_m$
be the  decomposition of $\Delta$ into left cosets, and  assume that
there are $\bbQ$-defined algebraic group morphisms 
$f_1, \ldots, f_m: \bG \rightarrow \bG$ 
so that 
\begin{equation} \label{exten}
 f_i(\gamma) = r_i \gamma r_i^{-1} \; , \; \, i= 1, \ldots, m, 
  \; \text{ for all } \gamma \in \Gamma \; . 
\end{equation} 
Then there exists a $\bbQ$-defined linear
algebraic group $\bI(\bG,\Gamma,\Delta)$, an injective  homomorphism 
$\psi: \Delta \longrightarrow \bI(\bG,\Gamma,\Delta)$ and 
a $\bbQ$-defined injective morphism of algebraic groups 
$\Psi: \bG \rightarrow  \bI(\bG,\Gamma,\Delta)$ 
which extends $\psi: \Gamma \rightarrow  \bI(\bG,\Gamma,\Delta)$ so that 
$\psi(\Delta) \leq  \bI(\bG,\Gamma,\Delta)_{\bbQ}$, $\bI(\bG,\Gamma,\Delta) \, = \, \Psi(\bG) \psi(\Delta)$ and $\Psi(\bG) \cap \psi(\Delta)= \psi(\Gamma)$. 
\end{lemma}

For more comments and details of the proof see \cite[Proposition 2.2]{GP1}.

\paragraph{The algebraic hull for $\Delta$}
We continue with our standing assumption that 
$\Delta$ is a virtually polycyclic group with $\Fitt(\Delta)$ 
torsion-free and $\C_{\Delta}(\Fitt(\Delta)) \leq \Fitt(\Delta)$.
Let $\Gamma \leq \Delta$ be a torsion-free polycyclic normal subgroup of
finite index. By Mostow's result there exists an algebraic hull 
$\bH_\Gamma$ for $\Gamma$, and we may assume that
$\Gamma \leq (\bH_\Gamma)_\bbQ$ is a Zariski-dense subgroup. 
Replacing $\Gamma$ with a
finite index subgroup, if necessary, we may 
also assume that $\bH_{\Gamma}$ is connected.  

\begin{proposition} \label{algreal} 
There exists a $\bbQ$-defined linear algebraic group $I(\bH_\Gamma,\Delta)$
which contains $\bH_\Gamma$, and an embedding
$\psi: \Delta \longrightarrow I(\bH_\Gamma,\Delta)_{\bbQ}$ 
which is the identity on $\Gamma$, such that
$I(\bH_\Gamma,\Delta)= \bH_\Gamma \psi(\Delta)$ 
and $\psi(\Delta) \cap \bH_\Gamma = \psi(\Gamma)$.
\end{proposition} 
\begin{proof} 
Let $\Delta= \Gamma r_1  \cup \, \cdots \, \cup \, \Gamma r_m$ be
a decomposition of $\Delta$ into left cosets. 
By the rigidity of the algebraic hull (Proposition 
\ref{rigidity1}), conjugation with $r_i$ on $\Gamma$ extends to
$\bbQ$-defined morphisms $f_i$ of $\bH$ 
which satisfy (\ref{exten}). 
The results follows then from Lemma \ref{induction}, 
putting $\bI(\bH_\Gamma,\Delta):= \bI(\bH_\Gamma,\Gamma,\Delta)$. 
\end{proof}

We need some more 
notations. Let $\bG$ be an algebraic 
group. We let $\bG^0$ denote its connected
component of identity. 
If $g$ is an element of $\bG$ then 
$g = g_u g_s$ denotes the Jordan-decomposition
of $g$ (i.e., $g_u$ is unipotent, $g_s$ is semisimple 
and $[g_s, g_u]=1$). If $M$ is a subset 
then let $M_u = \{ g_u \mid g \in M \}$,
$M_s=  \{ g_s \mid g \in M \}$. 
If $G$ is a subgroup then $\ur(G)$ denotes the unipotent 
radical of $G$, i.e., the maximal normal subgroup of $G$ which
consists of unipotent elements. 
We will use the following facts (compare, \cite[\S 10]{Borel}): 
If $G$ is a nilpotent
subgroup then $G_u$ and $G_s$ are subgroups of $\bG$, and 
$G \leq G_u \times G_s$. If $\bG^0$ is 
solvable then $\bG_u = \ur(\bG)$. 
In particular, for any subgroup $G$ of $\bG$, 
$\ur(G) = G \cap G_u$. \newline  

To construct the algebraic hull for $\Delta$
we have to further refine Proposition \ref{algreal}. 
\medskip 
 
\begin{prf}{Proof of Theorem \ref{alghull}} 
Let $\bU$ denote the unipotent radical of
$\bH_{\Gamma}$. By Proposition \ref{algreal}, we may assume that
$\Delta \leq \bG_{\bbQ}$ is a Zariski-dense subgroup of 
a $\bbQ$-defined linear algebraic group $\bG$ so that 
$\bG =\bH_\Gamma \Delta$, $\ur(\bG)=\bU$ and 
$\Delta \cap \bH_\Gamma = \Gamma$. 
Since $\bH_\Gamma$ is an algebraic hull,
$\Fitt(\Gamma) \leq \bU$, see Proposition \ref{Fittu}. 
Since $\Fitt(\Gamma)$ is a subgroup
of finite index in $\Fitt(\Delta)$ the group 
$\mu =\{ \gamma_s \mid \gamma \in \Fitt(\Delta) \} \leq \bG$ is finite.
Since $\Fitt(\Delta)$ is a normal subgroup of $\Delta$, 
$\mu$ is normalized by $\Gamma$. Hence, the centralizer of $\mu$ in $\bG$ 
contains a finite index subgroup of $\Gamma$.
Since the centralizer of $\mu$ is a Zariski-closed 
subgroup of $\bG$ it contains $(\bH_\Gamma)^0 = \bH_\Gamma$. 
In particular, $\mu$ centralizes $\Gamma$.  
We consider now the homomorphism $\psi_u: \Fitt(\Delta) 
\rightarrow \bU_{\bbQ}$    
which is given by $\gamma \mapsto \gamma_u$. 
The kernel of $\psi_u$ is contained in the finite group $\mu$.
Since $\Fitt(\Delta)$
is torsion-free, $\psi_u$ is injective.
Assigning $\psi: \delta \gamma \mapsto \psi_u(\delta) \gamma$ 
defines an injective homomorphism 
$\psi: \Fitt(\Delta) \Gamma \rightarrow (\bH_\Gamma)_\bbQ$. 
(To see that $\psi$ is injective suppose that $1 =\psi(\delta \gamma)$, for 
$\delta \in \Fitt(\Delta)$, $\gamma \in \Gamma$. 
Then $\gamma = \psi_u(\delta)^{-1}$ is unipotent, i.e., 
$\gamma \in \ur(\Gamma) \leq \Fitt(\Gamma)$. 
Therefore $\delta\gamma \in \Fitt(\Delta)$, and 
$\psi(\delta \gamma) = \psi_u(\delta \gamma)$,
and hence $\delta \gamma=1$.)
Clearly, the homomorphism $\psi$ is the identity on $\Gamma$.    
Let us put $\Gamma^*= \psi(\Fitt(\Delta) \Gamma)$.   
Then $\bH_\Gamma$ is an algebraic hull for $\Gamma^*$. 
We consider now the extension $\Gamma^* \leq \Delta$. 
By Proposition \ref{algreal}, there exist an algebraic group 
$I^*(\bH_\Gamma,\Delta)$ and an embedding 
$\psi^*: \Delta \longrightarrow \bI^*(\bH_\Gamma,\Delta)_{\bbQ}$  
so that  $\psi^*(\Delta) \cap \bI^*(\bH_\Gamma,\Delta) = \Gamma^*$. 
By construction, the group $\bI^*(\bH_\Gamma,\Delta)$ satisfies
i) and iii) of Theorem \ref{alghull} with respect to $\psi^*$
and $\Delta$. Moreover, by our construction $\psi^*(\Fitt(\Delta))
\leq \bU = \ur(\bI^*(\bH_\Gamma,\Delta))$.  

We consider now the centralizer $\C_{I^*(\bH_\Gamma,\Delta)}(\bU)$ 
of $\bU$ in $I^*(\bH_\Gamma,\Delta)$. It is a $\bbQ$-defined algebraic
subgroup of the solvable by finite group $I^*(\bH_\Gamma,\Delta)$. 
Therefore $\C_{I^*(\bH_\Gamma,\Delta)}(\bU) = \Z(\bU) S$, where
$S$ is a $\bbQ$-defined subgroup which 
consist of semisimple elements.
We remark that $S = \C_{I^*(\bH,\Delta)}(\bU)_s$ is normal 
in $I^*(\bH,\Delta)$. Since $S \cap \psi^*(\Delta)$ 
centralizes $\Fitt(\Delta) \leq \bU$, the assumption that
$\Delta$ has a strong Fitting subgroup implies that 
$S \cap \psi^*(\Delta) = \{1\}$. Let $\pi: I^*(\bH,\Delta) 
\rightarrow \bH_\Delta = I^*(\bH,\Delta)/S$ be the projection
homomorphism. Since $\bH_\Delta$ is $\bbQ$-defined with
a strong unipotent radical, $\bH_\Delta$ with the embedding 
$\pi \psi^*: \Delta \rightarrow (\bH_\Delta)_\bbQ$ is an algebraic hull 
for $\Delta$. 
\end{prf}


\subsection{Properties of the algebraic hull}
Let $\Gamma$ be a virtually polycyclic group. We 
assume that $\Gamma$ admits an algebraic hull $\bH_{\Gamma}$.  
The next proposition implies Proposition \ref{rigidity1},
and Corollary \ref{rigidity2}.

\begin{proposition}  \label{pchullextension}
Let $k$ be a subfield
of\/ $\bbC$, and  
$\bG$ a $k$-defined linear algebraic group with a strong
unipotent radical.  Let $\rho: \Gamma \longrightarrow \bG$ be 
a homomorphism so that $\rho(\Gamma)$ is Zariski-dense in $\bG$.
Then $\rho$ extends uniquely to a morphism of algebraic groups
$\rho_{\bH_{\Gamma}}: \bH_{\Gamma} \longrightarrow \bG$.
If $\rho(\Gamma) \leq   \bG_k$ then $\rho_{\bH_{\Gamma}}$ is 
defined over $k$.   
\end{proposition}
\begin{proof} We will use the diagonal argument. 
Therefore we consider  the subgroup 
$$ D= \{(\gamma, \rho(\gamma)) \mid \gamma \in \Gamma \} \;
 \leq \; \bH_\Gamma \times \bG \, \; .$$
Let $\pi_1, \pi_2$ denote the
projection morphisms onto the factors of the product  $\bH_\Gamma \times \bG$.
Let $\bD$ be the Zariski-closure of $D$,  and $\bU = \ur(\bD)$ the unipotent
radical of $\bD$. The group $\bD$ is a solvable by finite 
linear algebraic group, and  $\bD$ is defined over $k$
if $\rho(\Gamma) \leq \bG_k$.  
Let $\alpha=  {\pi_1}|_{\bD}: \bD \rightarrow \bH_\Gamma$. 
Since $\Gamma$ is Zariski-dense in $\bH_\Gamma$, 
$\alpha$ is onto. In particular, $\alpha$ maps $\bU= \bD_u$ 
onto $\ur(\bH_\Gamma)$.
By \cite[Lemma 4.36]{Raghunathan} we have 
$\dim \bU \leq \rank \Gamma = \dim \ur(\bH_\Gamma)$, and hence 
$\dim \bU = \dim \ur(\bH_\Gamma)$.  In particular, it follows
that the restriction $\alpha:  \bU \rightarrow \ur(\bH_\Gamma)$ 
is an isomorphism. 
Thus the kernel of $\alpha$ consists only of semi-simple
elements. Let $x \in \ker \,  \alpha$.  Then $x$ centralizes $\bU$.
Since $\pi_2(\bU)=\ur(\bG)$, $\pi_2(x)$ 
centralizes $\ur(\bG)$. Since 
$\bG$ has a strong unipotent radical, $x$ is in the kernel
of $\pi_2$, hence $x =1$. It follows that the morphism $\alpha$ 
is an isomorphism of groups. It is also an isomorphism of
algebraic groups. If $\rho(\Gamma) \leq \bG_k$ then  
$\alpha$ is $k$-defined. 
One can also show that $\alpha^{-1}$
is $k$-defined. 
(Compare e.g.\ \cite[Lemma 2.3]{GP1}.) 
We put $\rho_{\bH_{\Gamma}} = \pi_2 \circ \alpha^{-1}$
to get the required unique extension. $\rho_{\bH_{\Gamma}}$ 
is $k$-defined if $\rho(\Gamma) \leq \bG_k$. 
\end{proof}

\begin{remark} 
The proposition shows that the condition that $\ac{\rho(\Gamma)}$ has
a strong unipotent radical forces the homomorphism $\rho$ to be well
behaved.  For example, $\rho$ must be unipotent on the Fitting subgroup 
of $\Gamma$. See Proposition \ref{Fittu} below. 
\end{remark}

We study some further properties of the algebraic hull. 
In particular, we characterize the abstract virtually polycyclic groups 
which admit an algebraic hull in the sense of Theorem \ref{alghull}. 
Let us assume that $\Gamma$ is a Zariski-dense subgroup of a 
linear algebraic group $\bH$ with a strong 
unipotent radical. 

\begin{proposition} \label{Fittu} 
We have $\Fitt(\Gamma) \leq \ur(\bH)$.
In particular, $\ur(\Gamma) = \Fitt(\Gamma)$ and 
$\Fitt(\Gamma)$ is torsion-free. 
\end{proposition} 
\begin{proof} Let $F$ be the maximal 
nilpotent normal subgroup of $\bH$. 
Clearly, $F= \bF$ is a Zariski-closed subgroup.
Therefore  $\ur(\bF) = \ur(\bH)$. 
Now since $\bF$ is nilpotent,  $\bF_s$ is 
a subgroup, $\bF_u = \ur(\bF)$ and 
$\bF = \bF_s \cdot\ur(\bF)$ is a
direct product of groups. 
Since 
$\bH$ has a strong unipotent radical 
$\bF_s$ must be trivial, and it follows that $\bF= \ur(\bH)$.
The Zariski-closure of $\Fitt(\Gamma)$ is a 
nilpotent normal subgroup of $\bH$
and therefore $\Fitt(\Gamma)$ is contained in $\bF$.
Hence $\Fitt(\Gamma) \leq \ur(\bH)$. 
\end{proof} 

Let $\bN = \overline{\Fitt(\Gamma)}$ be the 
Zariski-closure of $\Fitt(\Gamma)$ in $\bH$. We just
proved $\bN \leq \ur(\bH)$.  
\begin{proposition} \label{FittCentralizer} 
Let $\bX= \C_\bH(\bN)$ be the centralizer of $\bN$
in $\bH$. Let $\bX^0$ be its component of identity. 
Then $\bX^0$ is a 
nilpotent normal 
subgroup of\/ $\bH$, and $\bX^0 \leq \ur(\bH)$. Moreover,
$\C_\Gamma(\Fitt(\Gamma)) \leq \Fitt(\Gamma)$. 
\end{proposition} 
\begin{proof} Since $\bX^0$ is 
a connected solvable algebraic group,
$\bX^0 = \bU \cdot \bT$, where
$\bU$ is a connected unipotent group and $\bT$ is
a maximal torus in $\bX^0$. Let $\bX_1= [\bX^0,\bX^0]$
and define $\Gamma_0= \Gamma \cap \bH^0$. Then $\Gamma_0$ is 
a polycyclic normal subgroup of $\Gamma$, and
$\Fitt(\Gamma_0) = \Fitt(\Gamma)$. Since
$\Gamma_0 \leq \bH^0$, it follows that 
$[\Gamma_0, \Gamma_0] \leq \ur(\Gamma_0) \leq \Fitt(\Gamma_0)$.  
This implies  
$[\bX^0,\bX^0] \leq [\bH^0, \bH^0]= \overline{[\Gamma_0, \Gamma_0]} \leq \bN$. 
We deduce $[\bX^0, \bX_1]= [\bT, \bX_1] = \{ 1\}$. 
On the other hand,
by \cite[\S 10.6]{Borel} all maximal tori in $\bX^0$ are
conjugate by an element of $\bX_1$. Hence $\bT$ must
be an invariant subgroup of $\bX^0$. In particular, 
$\bT$ is a normal abelian subgroup of $\bH$. Therefore, 
by the proof of the previous proposition, 
$\bT$ is contained in $\ur(\bH)$. Since $\bT$ consists
of semisimple elements, $\bT = \{1\}$. Hence, 
$\bX^0 = \Z(\ur(\bH))$.  

Next put $X= \C_\Gamma(\Fitt(\Gamma))$
and $X_0 = X \cap \bX^0$. 
Then $X_0$ is of finite index in $X$, 
nilpotent,  and a normal subgroup of $\Gamma$. 
The latter implies that $X_0 \leq \Fitt(\Gamma)$. 
It follows that 
$X_0$ is centralized by $X$. Since 
$X_0$ is of finite index in $X$ the commutator
subgroup $[X,X]$ must be finite. 
Since $[X,X]$ is normal in $\Gamma$ it follows that 
$\ur(\bH_\Gamma)$ centralizes $[X,X]$. Since 
$\bH$ has a strong unipotent radical 
$[X,X]= \{1\}$. It follows that $X$ is 
an abelian normal subgroup of $\Gamma$, 
and hence $X \leq \Fitt(\Gamma)$.   
\end{proof} 

Recall that $\Gamma$ is said to have a {\em strong Fitting subgroup}\/ 
if $\Fitt(\Gamma)$ is torsion-free and contains 
its centralizer. We summarize: 

\begin{corollary} \label{strongFitting}
Let $\Gamma$ be virtually polycyclic.
Then $\Gamma$ admits an algebraic hull $\bH$ if and
only if\/ $\Gamma$ has a strong Fitting subgroup.
\end{corollary}
We remark that this condition holds if and only if 
every finite normal subgroup of $\Gamma$ is trivial.
In fact, let us assume that the maximal normal 
finite subgroup of $\Gamma$ is trivial. Then 
$\Fitt(\Gamma)$ is torsion-free since its 
elements of finite order form a finite 
normal subgroup of $\Gamma$. Now put $X= \Z_\Gamma(\Fitt(\Gamma))$,
and let $X_0$ be a polycyclic normal subgroup of finite
index in $X$ which is nilpotent-by-abelian. 
From $[X_0, X_0] \leq \Fitt(X) \leq \Fitt(\Gamma)$  
we deduce that $X_0 \leq \Fitt(\Gamma)$. Therefore
$[X,X]$ must be finite, and it follows from our
assumption that $[X,X]= \{1 \}$. Hence, $X \leq \Fitt( \Gamma)$.



\end{appendix}


\begin{thebibliography}{99}

\bibitem[1]{Auslander}  L.\ Auslander, 
{An exposition of the structure of solvmanifolds. I. Algebraic theory}, 
 Bull.\  Amer.\  Math.\  Soc.\ {79} (1973), no.\ 2, 262-285.

\bibitem[2]{AJ} L.\ Auslander, F.E.A.\ Johnson,
{On a conjecture of C. T. C. Wall\/},
J.\ London Math.\ Soc.\ (2) {14} (1976), no.\ 2, 331-332.

\bibitem[3]{AusTol}  L.\ Auslander, R. Tolimieri,
 {On a conjecture of G. D. Mostow and the structure of solvmanifolds},
Bull.\ Amer.\ Math.\ Soc.\ {75}  (1969),  1330-1333.


\bibitem[4]{Baues1} O.\ Baues, {
Finite extensions and unipotent shadows of affine crystallographic groups}, 
C.\ R.\ Acad.\ Sci.\ Paris, Ser.\ I {335} (2002), 785-788.

\bibitem[5]{BenDek} Y.\ Benoist, K.\ Dekimpe,
{The uniqueness of polynomial crystallographic actions},
Math.\ Ann.\ {322} (2002),  no.\ 3, 563-571.

\bibitem[6]{Bieberbach} L.\ Bieberbach, {\"Uber die Bewegungsgruppen
der Euklidischen R\"aume II}, Math.\ Ann.\ {72} (1912), 400-412.

\bibitem[7]{Borel} A.\ Borel, {Linear algebraic groups},  
Second edition,  Graduate Texts in Mathematics {126},  
Springer-Verlag, (1991).

\bibitem[8]{BS} 
A.\ Borel,  J.-P.\  Serre, 
{Th\' eor\` emes de finitude en cohomologie galoisienne},
Comment.\ Math.\ Helv.\ { 39} (1964),  111-164.

\bibitem[9]{Browder} W.\ Browder,  {
On the action of $\Theta \sp{n}\,(\partial \pi )$},  
Differential and Combinatorial Topology 
(A Symposium in Honor of Marston Morse), 23-36, 
Princeton Univ.\  Press, Princeton, N.J., (1965).

\bibitem[10]{Cobb}  R.J.\ Cobb, {Infra-solvmanifolds of dimension four},
Bull.\ Austral.\ Math.\ Soc.\  {62} (2000),  
no.\ 2, 347-349.

\bibitem[11]{DenSing} Ch.\ Deninger, W.\ Singhof,
{On the cohomology of nilpotent Lie algebras},
Bull.\ Soc.\ Math.\ France {116} (1988), no.\ 1, 3-14.

\bibitem[12] {FJ1} F.T.\ Farrell, L.E.\ Jones, 
{ The surgery $L$-groups of poly-(finite or cyclic) groups},
Invent.\ Math.\ {91} (1988), no.\ 3, 559-586.

\bibitem[13]{FJ2} F.T.\ Farrell, L.E.\ Jones, 
{Compact infrasolvmanifolds are smoothly rigid},
Geometry from the Pacific Rim (Singapore, 1994),  85-97, 
de Gruyter, Berlin, (1997).

\bibitem[14]{FriedGoldman} D.\ Fried,  W.M.\ Goldman,
 {Three-dimensional affine crystallographic groups},  Adv.\ in Math.\ {47}
(1983), no.\ 1, 1-49.

\bibitem[15]{Goldman} W.M.\ Goldman, 
{On the polynomial cohomology of affine manifolds},
Invent.\ Math.\  {65}  (1981/82), no. 3, 453-457. 

\bibitem[16]{Hochschild} G.\ Hochschild, 
{Cohomology of algebraic linear groups},
Illinois J.\ Math.\  {5}  (1961),  492-519. 

\bibitem[17]{GP1}  F.\ Grunewald, V.\ Platonov,
{Solvable arithmetic groups and arithmeticity problems},  
Internat.\ J.\ Math.\  {10}   (1999),  no.\ 3, 327-366.

\bibitem[18]{GS}  F.\ Grunewald, D.\ Segal, 
{On affine crystallographic groups},
J.\  Differential Geom.\ { 40} (1994), no.\ 3, 563-594.

\bibitem[19]{IgDek} K.\ Dekimpe, P.\ Igodt,
{Polycyclic-by-finite groups admit a bounded-degree polynomial structure\/}, 
Invent.\  Math {129} (1997),  no.\ 1, 121-140.

\bibitem[20]{KS} R.C.\ Kirby, L.C.\ Siebenmann, 
{Foundational essays on topological manifolds, smoothings, and triangulations},
Annals of Mathematics Studies, No.\ 88,
Princeton University Press, Princeton, N.J., University of Tokyo Press, Tokyo, 1977.

\bibitem[21]{Lee-Raymond}  K.B.\ Lee, F.\ Raymond,
{Rigidity of almost crystallographic groups}, 
Combinatorial methods in topology and algebraic geometry (Rochester, N.Y., 1982),  73-78, 
Contemp.\ Math., {44}, Amer.\ Math.\ Soc., Providence, RI, 1985. 

\bibitem[22]{Lee-Raymond2}  K.B.\ Lee, F.\ Raymond,
{Maximal torus actions on solvmanifolds and double coset spaces},
Internat.\ J.\ Math.\ {2} (1991), no.\ 1, 67-76.

\bibitem[23]{Malcev} A.I.\ Malcev, {On a class of homogeneous spaces},
Amer.\ Math.\ Soc.\ Translation {39} (1951). 

\bibitem[24]{Mostow1} G.D.\ Mostow,
{Factor spaces of solvable groups\/},  
Ann.\ of Math.\ (2) {60} (1954),  1-27.

\bibitem[25]{Mostow3}  G.D.\ Mostow, 
{On the fundamental group of a homogeneous space}, 
Ann.\ of Math.\ (2) {66} (1957), 249-255.

\bibitem[26]{Mostow2}  G.D.\ Mostow, 
{Representative functions on discrete groups and solvable arithmetic subgroups}, 
Amer.\ J.\ Math.\  {92} (1970), 1-32. 

\bibitem[27]{Mostow5}  G.D.\ Mostow, 
Some applications of representative functions to solvmanifolds,
Amer.\ J.\ Math.\ {93} (1971) 11-32.

\bibitem[28]{Mostow4}  G.D.\ Mostow, 
{Strong rigidity of locally symmetric spaces},
Annals of Mathematics Studies, No.\ 78,
Princeton University Press, Princeton, N.J., 1973.  

\bibitem[29]{Nomizu}  K.\ Nomizu, {
On the cohomology of compact homogeneous spaces of nilpotent Lie groups},
Ann.\ of Math.\ (2)  {59} (1954), 531-538.

\bibitem[30]{Raghunathan} M.S.\ Raghunathan, 
{Discrete subgroups of Lie groups},  
Ergebnisse der Mathematik und ihrer Grenzgebiete {68},
Springer-Verlag, 1972. 

\bibitem[31]{Segal} D.\ Segal, {Polyclic groups}, Cambridge Univ.\ Press,
London, 1983.

\bibitem[32]{Tuschmann} W.\ Tuschmann, 
{Collapsing, solvmanifolds and infrahomogeneous spaces},
Differential Geom.\ Appl.\  {7} (1997),  no.\ 3, 251-264. 

\bibitem[33]{Wilking} B.\ Wilking, 
{Rigidity of group actions on solvable Lie groups},
Math.\ Ann.\  {317}  (2000),  no.\ 2, 195-237.


\end{thebibliography}
\end{document}